\documentclass{article}
\usepackage{latexsym}
\usepackage{epsfig}
\usepackage{amsmath}
\usepackage{amssymb}
\usepackage{amsthm}
\newtheorem{theorem}{Theorem}[section]
\newtheorem{corollary}[theorem]{Corollary}

\newtheorem{prop}[theorem]{Proposition}
\newtheorem{lemma}[theorem]{Lemma}

\theoremstyle{remark}
\newtheorem{remark}[theorem]{Remark}

\theoremstyle{definition}
\newtheorem{definition}[theorem]{Definition}

\DeclareMathOperator\Aut{Aut}

\DeclareMathOperator\diag{diag}
\DeclareMathOperator\id{id}
\DeclareMathOperator\rk{rk}

\DeclareMathOperator\sign{sgn}

\begin{document}

\title{Semidefinite descriptions of separable matrix cones}

\author{Roland Hildebrand \thanks{%
LJK, Tour IRMA, 51 rue des
Math\'ematiques, 38400 St.\ Martin d'H\`eres, France ({\tt
roland.hildebrand@imag.fr}). This paper presents research results
of the Belgian Programme on Interuniversity Poles of Attraction,
Phase V, initiated by the Belgian State, Prime Minister's Office
for Science, Technology and Culture; and of the Action Concert\'ee
Incitative "Masses de donn\'ees" of CNRS, France. The scientific
responsibility rests with its author.}}

\maketitle

\begin{abstract}
Let $K \subset E$, $K' \subset E'$ be convex cones residing in
finite-dimensional real vector spaces. An element $y$ in the tensor
product
$E \otimes E'$ is $K \otimes K'$-separable if it can be represented as
finite sum $y = \sum_l x_l \otimes x'_l$, where $x_l \in K$ and $x_l' \in
K'$ for all $l$.
Let ${\cal S}(n)$, ${\cal H}(n)$, ${\cal Q}(n)$ be the spaces of $n
\times n$ real symmetric, complex hermitian and quaternionic hermitian
matrices, respectively. Let further $S_+(n)$, $H_+(n)$, $Q_+(n)$ be the
cones of positive semidefinite matrices in these
spaces. If a matrix $A \in {\cal H}(mn) = {\cal H}(m) \otimes {\cal
H}(n)$ is $H_+(m) \otimes H_+(n)$-separable, then it fulfills also the
so-called PPT condition, i.e.\ it is positive semidefinite and has a
positive semidefinite partial transpose. 
The same implication holds for matrices in the spaces ${\cal S}(m) \otimes {\cal S}(n)$,
${\cal H}(m) \otimes {\cal S}(n)$, and for $m \leq 2$ in the space ${\cal Q}(m) \otimes {\cal S}(n)$.
We provide a complete enumeration of all pairs $(n,m)$ when the inverse implication is also true
for each of the above spaces, i.e.\ the PPT condition is sufficient for separability. 
We also show that a matrix in ${\cal Q}(n) \otimes {\cal S}(2)$ is
$Q_+(n) \otimes S_+(2)$-separable if and only if it is positive semidefinite.
\end{abstract}



\section{Introduction}

Let $K,K'$ be regular convex cones (closed convex cones, containing no
lines, with non-empty interior),
residing in finite-dimensional real vector spaces $E,E'$. Then an
element $w \in E \otimes E'$ of the tensor product space is called {\sl
$K \otimes K'$-separable}
(or just {\sl separable}, if it is clear which cones $K,K'$ are meant),
if it can be represented as a convex combination of product elements $v
\otimes v'$, where $v \in K$, $v' \in K'$.
It is not hard to show that the set of separable elements is itself a
regular convex cone. This cone is called the $K \otimes K'$-separable cone.

The notion of separability is intimately linked with the notion of
positive maps \cite{Choi74},\cite{Stormer63}. Cones of
positive maps appear frequently in applications
\cite{Terpstra},\cite{Horodeckis96} and are dual to separable cones
\cite{Gurvits0302102}.
Separability itself plays an increasingly important role in quantum
information theory \cite{Peres96}.

A particularly important case in optimization and Mathematical
Programming is when the cones $K,K'$ are standard self-scaled cones such
as Lorentz cones
or cones of positive semidefinite (PSD) matrices \cite{KoecherBraun}.
Let ${\cal S}(n)$, ${\cal H}(n)$, ${\cal Q}(n)$ be the spaces of $n
\times n$ real symmetric, complex hermitian and quaternionic hermitian
matrices, respectively. Let further $S_+(n)$, $H_+(n)$, $Q_+(n)$ be the
cones of PSD matrices in these spaces.
If the elements of a pair of matrix spaces commute, then the tensor product of these 
spaces can be represented by the Kronecker product space and is itself a subset of
such a matrix space. If such a product matrix is separable with respect to two PSD
matrix cones, then it is necessarily PSD itself.
In the case of $H_+(m) \otimes H_+(n)$-separability, where $n,m \in
{\mathbb N}_+$, there exists another simple necessary condition for
separability, the so-called PPT condition \cite{Peres96}.
A matrix in ${\cal H}(mn) = {\cal H}(m) \otimes {\cal H}(n)$ fulfills
the PPT condition if it is positive semidefinite and has a
positive semidefinite partial transpose.

In the spaces ${\cal S}(m) \otimes {\cal S}(n)$ and
${\cal H}(m) \otimes {\cal S}(n)$ the PPT condition reduces just to positivity, i.e.\ inclusion in the cone
$S_+(mn)$ or $H_+(mn)$, respectively. This is because the positivity property of real symmetric
or complex hermitian matrices is preserved under transposition.
However, the cone $Q_+(n)$ is invariant under transposition only for $n \leq 2$.
Therefore for matrices in ${\cal Q}(m) \otimes {\cal S}(n)$ the PPT condition is stronger than just positivity.
Moreover, it follows that the PPT condition is necessary for $Q_+(m) \otimes S_+(n)$-separability only for $m \leq 2$, while
positivity is necessary for arbitrary $(n,m)$.

The importance of the PPT condition and the positivity condition is based on the fact that these conditions
are semidefinite representable (i.e.\ in the form of linear matrix inequalities) and hence easily verifiable
algorithmically, in contrast to separability. It is then important to know in which cases these conditions are actually 
equivalent to separability, rather than only necessary. In these cases one then obtains semidefinite descriptions of the 
corresponding separable cones.

\medskip

The theorem of Woronowicz-Peres
states that in the case $m = 2$, $n = 3$ the PPT condition is not only
necessary, but also sufficient
for $H_+(m) \otimes H_+(n)$-separability
\cite{Woronowicz},\cite{Peres96}. However, there exist
matrices in ${\cal H}(2) \otimes {\cal H}(4)$ which fulfill the PPT condition, but are not separable \cite{Woronowicz}.
Similarly, Terpstra \cite{Terpstra} has shown that positivity is sufficient for $S_+(m) \otimes S_+(n)$-separability
for $\min(n,m) \geq 2$, but not for $n = m = 3$.
In \cite{Terpstra} this was formulated in the equivalent form of sums of squares representability of biquadratic forms.

One can then conclude that the positivity condition is equivalent to $S_+(m) \otimes S_+(n)$-separability
if and only if $\min(n,m) \leq 2$ and that the PPT condition is equivalent to $H_+(m) \otimes H_+(n)$-separability
if and only if $\min(n,m) = 1$ or $m+n \leq 5$.

In this contribution we provide a similar classification for the spaces ${\cal H}(m) \otimes {\cal S}(n)$ and ${\cal Q}(m) \otimes {\cal S}(n)$.
We show that positivity is equivalent to $H_+(m) \otimes S_+(n)$-separability if and only if $m = 1$ or $n \leq 2$ or $m+n \leq 5$,
and that the PPT condition is equivalent to $Q_+(2) \otimes S_+(n)$-separability if and only if $n \leq 3$.
Further, we show that for $m \geq 3$ positivity is equivalent to $Q_+(m) \otimes S_+(n)$-separability if and only if $n \leq 2$.
In addition, we enumerate all pairs $(n,m)$ for which the positivity property in ${\cal Q}(m) \otimes {\cal S}(n)$ is preserved by
the operation of matrix transposition, namely, the cases $m=1$, $n$ arbitrary, and the cases $m=2$, $n \leq 2$.
This involves mainly the following new and nontrivial results.

First, we show that a matrix in ${\cal Q}(2) \otimes {\cal S}(3)$ is
$Q_+(2) \otimes S_+(3)$-separable if and only if it fulfills the PPT
condition.
Second, we provide an example of a matrix in ${\cal H}(2) \otimes {\cal
S}(4)$ which fulfills the PPT condition, but is not $H_+(2) \otimes
S_+(4)$-separable, thus sharpening the counterexample
provided in \cite{Woronowicz} for the $H_+(2) \otimes H_+(4)$ case. 
Third, we show that if a matrix in ${\cal Q}(n) \otimes {\cal S}(2)$
is positive semidefinite, then it is $Q_+(n) \otimes S_+(2)$-separable.
In addition, we provide examples of matrices in ${\cal Q}(3)$ and ${\cal Q}(2) \otimes {\cal S}(3)$
which are PSD but whose transpose is not PSD.

\medskip

The remainder of the paper is structured as follows.
In the next section we provide exact definitions of separability and of
the PPT condition and consider some of their basic properties.
In Section 3 we consider low-dimensional cases and
relations between the cones we deal with. In the next two sections we prove
the sufficiency of the PPT condition for separability in
in the space ${\cal Q}(2) \otimes {\cal S}(3)$.
In section 5 we also provide an example of a PSD matrix in ${\cal Q}(2) \otimes {\cal S}(3)$
whose transpose is not PSD.
In Section 6 we provide a counterexample
against sufficiency of the PPT condition for $H_+(2) \otimes
S_+(4)$-separability.
In Section 7 we prove the equivalence of positivity and separability in the space
${\cal Q}(n) \otimes {\cal S}(2)$.
Finally we summarize our results in the last section.
In the appendix we list facts about quaternions and
quaternionic matrices which we use for the proof of the main
results of the paper. There we provide also an example of a matrix in $Q_+(3)$
whose transpose is not PSD.

\section{Definitions and preliminaries}

In this section we introduce the cones we deal with and provide
definitions of separability and the partial transpose. For basic
information related to quaternions and quaternionic matrices we refer
the reader to the appendix.
Throughout the paper, $i,j,k$ denote the imaginary units of the
quaternions and the overbar $\bar{\cdot}$ the complex or quaternion
conjugate. 

Let further $e_0,\dots,e_{m-1}$ be the canonical basis vectors of ${\mathbb R}^m$.
By $\id_E$ denote the identity operator on the space
$E$, by $I_n$ the $n \times n$ identity matrix, by $0_{n \times m}$ a zero matrix of size $n \times m$, and by $0_n$ a zero matrix
of size $n \times n$. Let further $\diag(A,B)$ denote a block-diagonal matrix with blocks $A$ and $B$.
For a matrix $A$ with real, complex or quaternionic entries,
$A^*$ will denote the transpose, complex conjugate transpose or quaternionic conjugate transpose of $A$, respectively,
and $\rk A$ the rank of $A$.
Further we denote by $GL_n(R)$ the set of invertible matrices of size $n \times n$ with entries in the ring $R$.

\medskip

We now introduce several convex cones we deal with.

Let $L_n$ be the $n$-dimensional standard Lorentz cone, or second order
cone,
\[ L_n = \left\{ (x_0,\dots,x_{n-1})^T \in \mathbb R^n \,|\, x_0 \geq
\sqrt{x_1^2 + \dots + x_{n-1}^2} \right\}.
\]
Let ${\cal S}(n)$ be the space of real symmetric $n \times n$ matrices and
$S_+(n)$ the cone of positive semidefinite (PSD) matrices in ${\cal S}(n)$;
${\cal H}(n)$ the space of complex hermitian $n \times n$
matrices and
$H_+(n)$ the cone of PSD matrices in ${\cal H}(n)$;
${\cal Q}(n)$ the space of quaternionic hermitian $n \times n$
matrices and
$Q_+(n)$ the cone of PSD matrices in ${\cal Q}(n)$;
${\cal Q}^k(n)$ the space of quaternionic hermitian $n \times n$
matrices with zero $k$-component and
$Q_+^k(n)$ the intersection of the PSD cone $Q_+(n)$ with the space
${\cal Q}^k(n)$.
All these cones are {\sl regular}, i.e.\ closed, containing no line, and with nonempty interior.

\medskip

Let $E,E'$ be real vector spaces and $K \subset E,K' \subset E'$
regular convex cones in
these spaces.

{\definition {\upshape \cite{Gurvits0302102}}
An element of the tensor product $E\otimes E'$ is called {\sl $K\otimes
K'$-separable}
if it can be written as a sum $\sum_{l=1}^N p_l x_l\otimes x'_l$, where
$N \in {\mathbb N}$ and $p_l > 0$,
$x_l \in K, x'_l \in K'$ for all $l = 1,\dots,N$. }

The set of $K\otimes K'$-separable elements forms a regular convex cone
in $E \otimes E'$.

\smallskip

We are interested in the case when the spaces $E,E'$ are spaces of hermitian matrices,
and the cones $K,K'$ are the corresponding cones of PSD matrices in these spaces.
If the elements in the factor spaces commute, then we can represent tensor products 
of matrices by Kronecker products, which will again be hermitian matrices. If the elements
do not commute, then the Kronecker products will in general not be hermitian.
Since we are interested in describing the separable cone by a matrix inequality,
we do not consider this latter case in this contribution.

In particular, for $m,n \in {\mathbb N}_+$ we consider the space ${\cal S}(m) \otimes {\cal S}(n)$
as a subspace of ${\cal S}(mn)$, the spaces ${\cal H}(m) \otimes {\cal S}(n)$ and
${\cal H}(m) \otimes {\cal H}(n)$ as subspaces of ${\cal H}(mn)$ and the space
${\cal Q}(m) \otimes {\cal S}(n)$ as subspace of ${\cal Q}(mn)$. As can easily be seen,
exchanging the factors in the tensor product is equivalent to applying a certain permutation
of rows and columns to the corresponding Kronecker product. Hence
exchanging the factors in the tensor product leads to a canonically isomorphic space,
and we can restrict our consideration to the cases listed above.

{\definition An automorphism $A$ of $E$ is called an {\sl automorphism
of $K$} if $A[K] = K$. }

\smallskip

The automorphisms of a cone $K$ form a group, which will be called
$\Aut(K)$. For regular convex cones $K,K'$ let $g \in \Aut(K),g' \in
\Aut(K')$ be elements of their automorphism groups. Since $g,g'$ are
linear automorphisms of the underlying spaces $E,E'$, we can consider
their tensor product $g \otimes g'$, which will be a linear automorphism
of the space $E \otimes E'$.
The following assertion is trivial but nevertheless very useful.

{\prop \label{prodautgroup}
Let ${\cal G}_{K \otimes K'}$ be the set of maps $\{ g \otimes g':E
\otimes E' \to E \otimes E' \,|\, g \in \Aut(K),g' \in \Aut(K') \}$.
Then ${\cal G}_{K \otimes K'}$ is a group. It is canonically isomorphic to the
factor group
$[\Aut(K) \times \Aut(K')]/\{(\alpha \id_V,\alpha^{-1}\id_{V'})\,|\,\alpha \in {\mathbb R}_+ \}$. 
The cone of $K \otimes K'$-separable elements
is invariant with respect to the action of ${\cal G}_{K \otimes K'}$,
and ${\cal G}_{K \otimes K'}$ is a subgroup of the automorphism group of this cone. \qed}

{\definition A {\sl face} $F$ of a convex cone is a subset of $K$ with
the following property. If $x,y \in K$ and $\frac{x+y}{2} \in F$, then
$x,y \in F$. For $x \in K$, {\sl the face} of $x$ in $K$
is the minimal face of $K$ containing $x$. }

It is not hard to see that a face $F$ is the face of a point $x$ if and only if $x$ is contained in the relative interior of $F$. If $x \in K$ and $y \in K'$, then the face of $(x,y)$ in $K \times K' \subset E \times E'$
is $F_x \times F_y$, where $F_x$ is the face of $x$ in $K$ and $F_y$ is the face of $y$ in $K'$.

{\definition An {\sl extreme ray} of a regular convex cone $K$ is a 1-dimensional
face. Any non-zero point on an extreme ray is a {\sl generator} of that
extreme ray. }

A convex cone is the convex hull of its extreme rays.

{\lemma \label{dimcount} Let $K \subset {\mathbb R}^N$ be a regular
convex cone and $L
\subset {\mathbb R}^N$ a linear subspace of dimension $n$. Let $K' = K
\cap L$ and let $x$ be the generator of an extreme ray of $K'$.
Then the face of $x$ in $K$ has at most dimension $N - n +1$. }

\begin{proof} Let $d$ be the dimension of the face $F$ of $x$ in
$K$. Then the intersection $F \cap L$ has dimension $d' \geq d + n
- N$ and $x$ is in the relative interior of $F \cap L \subset K'$. On the other hand, $d' = 1$, since $x$ generates an extreme ray of
$K'$. It follows that $d \leq N - n + 1$.
\end{proof}

Any face of $S_+(n)$, $H_+(n)$, or $Q_+(n)$ is isomorphic to the cone
$S_+(l)$, $H_+(l)$, or $Q_+(l)$, respectively,  for some $l \in \{0,\dots,n\}$, has
dimension $\frac{l(l+1)}{2}$, $l^2$, or
$l(2l-1)$, respectively, and consists of matrices of rank not exceeding $l$.
The proof of this assertion is similar for all three cases, for the quaternionic case we refer
the reader to Proposition \ref{quatface} in the appendix.

\medskip

Denote the spaces ${\cal S}(2)$, ${\cal H}(2)$, ${\cal Q}^k(2)$, ${\cal Q}(2)$ by
$E_3,E_4,E_5,E_6$, and the cones $S_+(2)$, $H_+(2)$, $Q_+^k(2)$, $Q_+(2)$ by
$K_3,K_4,K_5,K_6$, respectively. The index denotes the real dimension of
the corresponding space or cone.
Let $T_m: E_m \to E_m$, $m = 3,4,5,6$, be the matrix transposition, or equivalently the
complex or quaternion conjugation in $E_m$.
Note that $T_m$ is in $\Aut(K_m)$.

\smallskip

We have the inclusions $E_m \subset E_n$ and $K_m \subset K_n$ for $m
\leq n$. For $n \in {\mathbb N}_+$ we have
$E_3 \otimes {\cal S}(n) \subset {\cal S}(2n)$, $E_4 \otimes {\cal S}(n)
\subset {\cal H}(2n)$, $E_5 \otimes {\cal S}(n) \subset {\cal Q}^k(2n)$
and $E_6 \otimes {\cal S}(n) \subset {\cal Q}(2n)$. The space $E_m \otimes {\cal S}(n)$, $m = 3,\dots,6$,
consists of matrices composed of 4 symmetric $n \times n$ blocks.

{\definition Let $m,n \in {\mathbb N}_+$ and let $A$ be a $mn \times mn$ matrix, with real, complex or
quaternionic entries. Partition $A$ in $m \times m$ blocks $A_{\alpha\beta}$ ($\alpha,\beta = 1,\dots,m$)
of size $n \times n$.
Then the {\sl partial transpose} of $A$, denoted by $A^{\Gamma}$, will
be defined as the result of exchanging the off-diagonal blocks $A_{\alpha\beta}$ and $A_{\beta\alpha}$
for all $\alpha \not= \beta$.
\[ A = \left( \begin{array}{cccc} A_{11} & A_{12} & \cdots & A_{1m} \\ A_{21} & A_{22} & \cdots & A_{2m} \\
\vdots & \vdots & & \vdots \\ A_{m1} & A_{m2} & \cdots & A_{mm}
\end{array} \right), \qquad
A^{\Gamma} = \left( \begin{array}{cccc} A_{11} & A_{21} & \cdots & A_{m1} \\ A_{12} & A_{22} & \cdots & A_{m2} \\
\vdots & \vdots & & \vdots \\ A_{1m} & A_{2m} & \cdots & A_{mm}
\end{array} \right).
\]
}

Note that if $A$ is hermitian, i.e.\ $A = A^*$, then so is $A^{\Gamma}$.

{\definition Let $A$ be a hermitian $mn \times mn$ matrix, with real,
complex or quaternionic entries. Then $A$ is said to fulfill the {\sl
PPT condition} or to be a {\sl PPT matrix}
if both $A$ and $A^{\Gamma}$ are positive semidefinite. }

The set of PPT matrices in the spaces ${\cal S}(m) \otimes {\cal S}(n)$,
${\cal H}(m) \otimes {\cal S}(n)$, ${\cal Q}(m) \otimes {\cal S}(n)$, ${\cal H}(m) \otimes {\cal H}(n)$
for fixed $m,n \in {\mathbb N}_+$ is a regular convex cone.

It is well-known that for matrices in ${\cal H}(mn)$ being PPT is a necessary condition
for $H_+(m) \otimes H_+(n)$-separability \cite{Peres96}. 
We can generalize this result in the following way.

{\prop \label{pospptnec} Let $V_m$ be one of the matrix spaces ${\cal S}(m),{\cal H}(m),{\cal Q}(m)$ and $V_n$ one of the spaces 
${\cal S}(n),{\cal H}(n),{\cal Q}(n)$, such that the elements of $V_m$ and $V_n$ commute.
Let $K_+(m),K_+(n)$ be the corresponding positive matrix cones.
Then any $K_+(m) \otimes K_+(n)$-separable matrix in $V_m \otimes V_n$ is PSD.

Assume further that $K_+(m)$ is invariant with respect to transposition. Then
any $K_+(m) \otimes K_+(n)$-separable matrix in $V_m \otimes V_n$ is a PPT matrix. }

\begin{proof}
It is sufficient to prove the assertions for the extreme rays of the $K_+(m) \otimes K_+(n)$-separable cone.
Let the matrix $A$ generate such an extreme ray. Then $A$ can be written as Kronecker product
$A_m \otimes A_n$, where $A_m = xx^* \in K_+(m)$, $A_n = yy^* \in K_+(n)$ generate extreme rays
of the corresponding PSD cones, and $x,y$ are appropriate column vectors
of size $m,n$, respectively. Hence $A = (x \otimes y)(x \otimes y)^*$ is PSD, which proves the first part
of the proposition.

Let us prove the second part. We have $A^{\Gamma} = A_m^T \otimes A_n$. By assumption of the proposition
$A_m^T$ also generates an extreme ray of $K_+(m)$ and can hence be expressed as $A_m^T = x'{x'}^*$, where
$x'$ is an appropriate column vector of length $m$. Then $A^{\Gamma} = (x' \otimes y)(x' \otimes y)^*$ is also PSD
and $A$ is a PPT matrix.
\end{proof}

In particular, the PPT property is necessary for $S_+(m) \otimes S_+(n)$-, $H_+(m) \otimes S_+(n)$-, and 
$H_+(m) \otimes H_+(n)$-separability for any $m,n \in {\mathbb N}_+$, and for $Q_+(m) \otimes S_+(n)$-separability
for $m \leq 2$ and $n \geq 1$ arbitrary (cf.\ Corollary \ref{posQ2} in the appendix).

Note also that if $V_n = {\cal S}(n)$, then the operation of partial transposition is equivalent to full transposition.
If in addition, transposition preserves positivity of matrices in $V_m \otimes V_n$, then a matrix has the PPT property
if and only if it is PSD. In particular, this holds for the spaces ${\cal S}(m) \otimes {\cal S}(n)$ and ${\cal H}(m) \otimes {\cal S}(n)$.

\smallskip

Let now $3 \leq m \leq 6$ and $n \geq 1$.
Denote by $\Sigma_{m,n}$ the cone of $K_m \otimes S_+(n)$-separable matrices
and by $\Gamma_{m,n}$ the cone of PPT matrices in the space $E_m \otimes
{\cal S}(n)$. Observe that $K_m = \Sigma_{m,1} = \Gamma_{m,1}$ for all $m$.
Proposition \ref{pospptnec} yields the following result.

{\corollary For any $3 \leq m \leq 6$ and $n \in {\mathbb N}_+$, we have the
inclusion $\Sigma_{m,n} \subset \Gamma_{m,n}$. \qed }

Let ${\cal G}_{m,n}$ be the group $\{ g \otimes g' \,|\, g \in
\Aut(K_m),g' \in \Aut(S_+(n)) \}$. By Proposition \ref{prodautgroup} it is a subgroup of $\Aut(\Sigma_{m,n})$.
Note that the operator $T_m \otimes \id_{{\cal S}(n)}$ of partial
transposition is in ${\cal G}_{m,n}$ and amounts to complex or
quaternion conjugation.

Let us define isomorphisms ${\cal I}_m: {\mathbb R}^m \to E_m$, $3 \leq m \leq 6$.
\begin{eqnarray} \label{LorIso}
{\cal I}_3: && (x_0,x_1,x_2)^T \mapsto \left( \begin{array}{cc} x_0+x_1 & x_2 \\ x_2 &
x_0-x_1 \end{array} \right) \nonumber\\
{\cal I}_4: && (x_0,x_1,x_2,x_3)^T \mapsto \left( \begin{array}{cc} x_0+x_1 & x_2+ix_3
\\ x_2-ix_3 & x_0-x_1 \end{array} \right) \nonumber\\
{\cal I}_5: && (x_0,x_1,x_2,x_3,x_4)^T \mapsto \left( \begin{array}{cc} x_0+x_1 &
x_2+ix_3+jx_4 \\ x_2-ix_3-jx_4 & x_0-x_1 \end{array} \right) \nonumber\\
{\cal I}_6: && (x_0,x_1,x_2,x_3,x_4,x_5)^T \mapsto \left( \begin{array}{cc} x_0+x_1 &
x_2+ix_3+jx_4+kx_5 \\ x_2-ix_3-jx_4-kx_5 & x_0-x_1 \end{array} \right)
\end{eqnarray}

It is not hard to check the following result (cf.\ Corollary \ref{posQ2} in the appendix).

{\lemma $K_m = {\cal I}_m[L_m]$ for all $m = 3,4,5,6$. \qed}

As a consequence, the group $\Aut(K_m)$ is
isomorphic to the automorphism group of the Lorentz cone $L_m$.

By virtue of (\ref{LorIso}) the map ${\cal I}_m \otimes \id_{{\cal S}(n)}$ is an isomorphism between the spaces
${\mathbb R}^m \otimes {\cal S}(n)$ and $E_m \otimes {\cal S}(n)$. Hence we can represent any element of $E_m \otimes {\cal S}(n)$
in a unique way as an image $({\cal I}_m \otimes \id_{{\cal S}(n)})(\sum_{l=0}^{m-1} e_l \otimes B_l)$,
where $B_0,\dots,B_{m-1} \in {\cal S}(n)$.

{\definition For any element ${\bf B} = ({\cal I}_m \otimes \id_{{\cal S}(n)})(\sum_{l=0}^{m-1} e_l \otimes B_l) \in E_m \otimes {\cal S}(n)$,
the matrices $B_0,\dots,B_{m-1} \in {\cal S}(n)$ will be called the {\sl components} of ${\bf B}$. }

{\prop \label{symmetry} The cone $\Gamma_{m,n}$ is invariant with
respect to the action of ${\cal G}_{m,n}$, i.e.\ ${\cal G}_{m,n} \subset
\Aut(\Gamma_{m,n})$. }

Before proceeding to the proof of this proposition, we define
$\gamma_m$, $m = 3,\dots,6$ to be the set of all appropriate matrices
$S$ such that the mapping $A \mapsto SAS^*$ is an automorphism of the
space $E_m$.
More precisely, define $\gamma_3 = GL_2({\mathbb R}), \gamma_4 =
GL_2({\mathbb C}), \gamma_6 = GL_2({\mathbb H})$ and let $\gamma_5$ be
the set of all matrices $S \in GL_2({\mathbb H})$ such that
$SAS^* \in E_5$ whenever $A \in E_5$. The set $\gamma_5$ is a matrix
group, its Lie algebra given by all quaternionic $2 \times 2$ matrices
which: i) have a trace with zero $i$- and $j$-components, ii) the
$k$-components of the off-diagonal elements are zero, and iii) the
$k$-components of the diagonal elements are equal.

Observe that mappings of the form $A \mapsto SAS^*$ preserve the PSD
matrix cone $K_m$. Let then $H^{\gamma}_m: \gamma_m \to \Aut(K_m)$ be
the group homomorphism assigning to any
matrix $S \in \gamma_m$ the automorphism $H^{\gamma}_m(S): A \mapsto
SAS^*$. Denote the image of $H^{\gamma}_m$ by $G_m$. Hence $G_m$ is the
subgroup of all automorphisms of $E_m$ of the form
$A \mapsto SAS^*$. Let also $H^{{\mathbb R}}_n: GL_n({\mathbb R}) \to
\Aut(S_+(n))$ be the group homomorphism assigning to any matrix $S \in
GL_n({\mathbb R})$ the automorphism
$H^{{\mathbb R}}_n(S): A \mapsto SAS^T$. This group homomorphism is
surjective.

\begin{proof}[Proof of Proposition \ref{symmetry}]
Assume the notations of the proposition.

The group $\Aut(L_m)$ is a Lie group of dimension $\frac{m(m-1)}{2}+1$
and consists of two connected components, corresponding to automorphisms
with positive and negative determinants, respectively. Hence $\Aut(K_m)$
is generated by the connected component of its neutral element and
a single automorphism with negative determinant. By computing the rank
of the Lie algebra of $G_m$ one easily determines that the dimension of
$G_m$ equals that of $\Aut(L_m)$.
It follows \cite{Warner} that $G_m$ contains the connected component of
the neutral element of $\Aut(K_m)$. For $m = 3,5$ the group $G_m$
contains automorphisms with negative determinants and hence
equals $\Aut(K_m)$. For $m = 4,6$ $G_m$ is connected, because $\gamma_m$ is
connected. However, in this case the matrix transposition $T_m$ has
negative determinant as an automorphism of $E_m$.
Thus $\Aut(K_m)$ is generated by $G_m$ and $T_m$ for all $m = 3,\dots,6$.

Therefore ${\cal G}_{m,n}$ is generated by the following elements.
First, elements of the form $H^{\gamma}_m(S) \otimes \id_{{\cal S}(n)}$,
where $S \in \gamma_m$; second, the element
$T_m \otimes \id_{{\cal S}(n)}$; and third, elements of the form
$\id_{E_m} \otimes H^{{\mathbb R}}_n(S)$, where $S \in GL_n({\mathbb R})$.
Let us now consider the action of these generators on the cone
$\Gamma_{m,n}$. Let $K_+ \subset E_m \otimes {\cal S}(n)$ be the cone of
PSD matrices in $E_m \otimes {\cal S}(n)$.

Let $S \in \gamma_m$. Then $T_m^{-1} \circ H^{\gamma}_m(S) \circ T_m$ is
an element of $G_m$, and there exists a matrix $S' \in \gamma_m$ such
that $T_m^{-1} \circ H^{\gamma}_m(S) \circ T_m = H^{\gamma}_m(S')$.
The element $H^{\gamma}_m(S) \otimes \id_{{\cal S}(n)}$ acts on $E_m
\otimes {\cal S}(n)$ as $A \mapsto (S \otimes I_n) A (S \otimes I_n)^*$
and hence preserves the cone $K_+$.
Moreover, we have $(H^{\gamma}_m(S) \otimes \id_{{\cal S}(n)}) \circ
(T_m \otimes \id_{{\cal S}_n}) = (H^{\gamma}_m(S) \circ T_m) \otimes
\id_{{\cal S}_n} =
(T_m \circ H^{\gamma}_m(S')) \otimes \id_{{\cal S}_n} = (T_m \otimes
\id_{{\cal S}_n}) \circ (H^{\gamma}_m(S') \otimes \id_{{\cal S}(n)})$.
If we denote the set $\{A^{\Gamma} \,|\, A \in K_+\}$ by $K_+^{\Gamma}$,
then
we get $(H^{\gamma}_m(S) \otimes \id_{{\cal S}(n)})[K_+^{\Gamma}] =
(H^{\gamma}_m(S) \otimes \id_{{\cal S}(n)}) \circ (T_m \otimes
\id_{{\cal S}_n})[K_+] = (T_m \otimes \id_{{\cal S}_n}) \circ
(H^{\gamma}_m(S') \otimes
\id_{{\cal S}(n)})[K_+] = (T_m \otimes \id_{{\cal S}_n})[K_+] =
K_+^{\Gamma}$. Hence $H^{\gamma}_m(S) \otimes \id_{{\cal S}(n)}$
preserves also the cone $K_+^{\Gamma}$, and therefore also
$\Gamma_{m,n} = K_+ \cap K_+^{\Gamma}$.

The automorphism $(T_m \otimes \id_{{\cal S}_n})$ is the operator of
partial transposition on $E_m \otimes \id_{{\cal S}(n)}$ and hence
preserves $\Gamma_{m,n}$ by definition.

Let now $S \in GL_n({\mathbb R})$. Then $\id_{E_m} \otimes H^{{\mathbb
R}}_n(S)$ acts on $E_m \otimes {\cal S}(n)$ like
$A \mapsto (I_2 \otimes S)A(I_2 \otimes S)^T$ and hence preserves $K_+$.
Moreover, $\id_{E_m} \otimes H^{{\mathbb R}}_n(S)$ commutes with the
operator $T_m \otimes \id_{{\cal S}(n)}$ of partial transposition and hence
preserves also $K_+^{\Gamma}$. Therefore it preserves $\Gamma_{m,n}$.

Thus all generators of ${\cal G}_{m,n}$ preserve the cone
$\Gamma_{m,n}$, and so do all other elements.
\end{proof}

Proposition \ref{symmetry} allows us, when looking for elements in
$\Gamma_{m,n}\setminus \Sigma_{m,n}$,
to restrict the consideration to elements that are in some canonical
form with respect to the action
of the symmetry group ${\cal G}_{m,n}$, since the
inclusions in $\Sigma_{m,n}$ or $\Gamma_{m,n}$ hold or do not hold
for all elements of an orbit simultaneously.

\section{Relations between different cones and trivial cases}

The next three sections aim at proving that the cones
$\Sigma_{6,3}$ and $\Gamma_{6,3}$ are equal. 

{\prop \label{reduceimpl} Let $3 \leq m' \leq m \leq 6$ and $n \geq n' >
0$ and suppose that $\Gamma_{m,n} = \Sigma_{m,n}$. Then $\Gamma_{m',n'}
= \Sigma_{m',n'}$. }

\begin{proof} Assume the conditions of the proposition. It is sufficient
to prove the assertion for $n' = n-1,m' = m$ and $n' = n,m' = m-1$.

Let first $n' = n-1,m' = m$. Define the linear mapping $L_{n',n}: {\cal
S}(n') \to {\cal S}(n)$ as follows. For $S \in {\cal S}(n')$, let
$L_{n',n}(S)$ be the matrix which has its upper left $n' \times n'$
submatrix equal to $S$
and whose all other elements are zero. It is not hard to see that if a
matrix $A \in E_m \otimes {\cal S}(n)$ is in the image of $\id_{E_m} \otimes L_{n',n}$,
then $A \in \Sigma_{m,n}$ if and only if $A \in (\id_{E_m} \otimes L_{n',n})[\Sigma_{m,n'}]$
and $A \in \Gamma_{m,n}$ if and only if $A \in (\id_{E_m} \otimes L_{n',n})[\Gamma_{m,n'}]$.
Injectivity of $L_{n',n}$ now yields the desired result.

Let now $n' = n,m' = m-1$ and let $C \in \Gamma_{m',n}$. Then we have
also $C \in \Gamma_{m,n}$ and by assumption $C \in \Sigma_{m,n}$.
Therefore we can represent $C$ as a sum
$\sum_{l=1}^N A_l \otimes B_l$, where $A_l \in K_m$ and $B_l \in
S_+(n)$. Let now $\pi_m: {\mathbb R}^m \to {\mathbb R}^{m'}$ be the projection that
assigns to any vector $(x_0,\dots,x_{m-2},x_{m-1})^T \in {\mathbb R}^m$ the vector
$(x_0,\dots,x_{m-2})^T \in {\mathbb R}^{m'}$ and let 
$\pi_m^E = {\cal I}_{m'} \circ \pi_m \circ {\cal I}^{-1}_m: E_m \to E_{m'}$.
Since $\pi_m[L_m] = L_{m'}$, we have $\pi_m^E[K_m] = K_{m'}$.
Moreover, the restriction of $\pi_m^E$ to $E_{m'}$ is just $\id_{E_{m'}}$
and the restriction of $\pi_m^E \otimes {\cal S}(n)$ to $E_{m'} \otimes {\cal S}(n)$ is just $\id_{E_{m'} \otimes {\cal S}(n)}$.
Therefore $C = (\pi_m^E \otimes {\cal S}(n))(C) = \sum_{l=1}^N A_l' \otimes B_l$
with $A_l' = \pi_m^E(A_l)$. But $A_l'$ is an element of $K_{m'}$, because
$A_l \in K_m$. Therefore $C$ is
$K_{m'}\otimes S_+(n)$-separable and hence in $\Sigma_{m',n}$.
\end{proof}

As one should expect, the equality $\Gamma_{m,n} = \Sigma_{m,n}$ is thus
easier to prove for smaller $n,m$. 

By a similar reasoning we can prove the following results.

{\lemma \label{33conseq} Let $n,m \geq 3$. Then the positivity of a matrix in ${\cal H}(m) \otimes {\cal S}(n)$ is not sufficient
for $H_+(m) \otimes S_+(n)$-separability. The PPT condition in ${\cal Q}(m) \otimes {\cal S}(n)$ is not sufficient for
$Q_+(m) \otimes S_+(n)$-separability. }

\begin{proof}
Let $n,m \geq 3$. Then the cone of positive semidefinite matrices in the space ${\cal S}(m) \otimes {\cal S}(n)$ does not coincide
with the cone of $S_+(m) \otimes S_+(n)$-separable matrices. This is a consequence of the fact that not every nonnegative definite biquadratic form $F(x,y)$,
where $x \in {\mathbb R}^m$, $y \in {\mathbb R}^n$ are two vectors of variables, is representable as a sum of squares of bilinear forms \cite{Terpstra}.
Let then $C \in {\cal S}(m) \otimes {\cal S}(n)$ be a PSD matrix which is not $S_+(m) \otimes S_+(n)$-separable. Then $C$ is
also PSD if considered as an element of ${\cal H}(m) \otimes {\cal S}(n)$ and fulfills the PPT condition if considered as an element of
${\cal Q}(m) \otimes {\cal S}(n)$.

Now suppose that $C$ is $H_+(m) \otimes S_+(n)$-separable. Then there exists an integer $N \in {\mathbb N}$ and matrices $A_1,\dots,A_N \in H_+(m)$,
$B_1,\dots,B_N \in S_+(n)$ such that $C = \sum_{l=1}^N A_l \otimes B_l$. Since $C$ is real, we then have
$C = \sum_{l=1}^N Re(A_l) \otimes B_l$. But $Re(A_l) \in S_+(m)$ for all $l$. Therefore $C$ is $S_+(m) \otimes S_+(n)$-separable, contrary
to our assumption.

In the same way one shows that $Q_+(m) \otimes S_+(n)$-separability of $C$ leads to a contradiction.
\end{proof}

It now happens that an element of $\Sigma_{m,n}$ can be reduced to elements of
cones in lower dimensions. We formalize this in the following definitions.

{\definition \label{def3} Let $3 \leq m \leq 6$ and $n \geq 1$. Let ${\bf B} \in E_m \otimes {\cal S}(n)$ have components $B_0,\dots,B_{m-1} \in {\cal S}(n)$. 
If there exists a matrix $S \in GL_n({\mathbb R})$ and positive integers $n_1,n_2$ with $n_1+n_2 = n$ such that the matrices $SB_lS^T$ are
block-diagonal with blocks $B^1_l,B^2_l$ of sizes $n_1\times n_1,n_2\times n_2$ for all $l = 0,\dots,m-1$,
then we call the element ${\bf B}$ {\sl decomposable}. We call the elements ${\bf B}^1 = ({\cal I}_m \otimes {{\cal S}(n_1)})(\sum_{l=0}^{m-1} e_l \otimes B_l^1)$, 
${\bf B}^2 = ({\cal I}_m \otimes {{\cal S}(n_2)})(\sum_{l=0}^{m-1} e_l \otimes B_l^2)$ {\sl block components} of ${\bf B}$. }

An element of $E_m \otimes {\cal S}(n)$ may have several decompositions and its block components are not uniquely defined.

{\definition Let $4 \leq m \leq 6$ and $n \geq 1$. If an element ${\bf B} \in E_m \otimes {\cal S}(n)$ has linearly dependent components, then we call it {\sl reducible}. 
If there exists $A \in \Aut(K_m)$ such that ${\bf B}' = (A \otimes \id_{{\cal S}(n)})({\bf B}) \in E_{m-1} \otimes {\cal S}(n)$, then we call
${\bf B}'$ a {\sl reduction} of ${\bf B}$. }

{\lemma Let ${\bf B} \in \Gamma_{m,n}$ be reducible. Then ${\bf B}$ has a reduction. }

\begin{proof}
Assume the conditions of the lemma and let $B_0,\dots,B_{m-1}$ be the components of ${\bf B}$. 
Then there exists a nonzero vector $\lambda = (\lambda_0,\dots,\lambda_{m-1})^T \in {\mathbb R}^m$ such that
$\sum_{l=0}^{m-1} \lambda_l B_l = 0$. 

If ${\bf B} = 0$, then ${\bf B} \in E_{m-1} \otimes {\cal S}(n)$ and ${\bf B}$ is a reduction of itself.
Hence let us assume that ${\bf B} \not= 0$. Since ${\bf B} \in \Gamma_{m,n}$, the diagonal elements of
the components $B_0,\dots,B_{m-1}$ cannot all be zero simultaneously. Without restriction of generality
let us assume that the $(1,1)$-elements of $B_0,\dots,B_{m-1}$ form a nonzero vector
$b = (b_0,\dots,b_{m-1})^T \in {\mathbb R}^m$. 
Define $b' = (b_1,\dots,b_{m-1})^T, \lambda' = (\lambda_1,\dots,\lambda_{m-1})^T \in {\mathbb R}^{m-1}$.

The element ${\bf b} = {\cal I}_m(b) \in E_m$ is a
principal $2 \times 2$ submatrix of ${\bf B}$. Since ${\bf B} \in \Gamma_{m,n}$, we have ${\bf b} \succeq 0$ and hence ${\bf
b} \in K_m$, $b \in L_m$, and $0 \not= b_0 \geq ||b'||$.
Moreover, $\langle b, \lambda \rangle = b_0\lambda_0 + \langle b',\lambda' \rangle = 0$.
It follows that $b_0|\lambda_0| \leq ||b'||\,||\lambda'|| \leq b_0 ||\lambda'||$ and $|\lambda_0| \leq ||\lambda'||$.
We distinguish two cases.

i) $|\lambda_0| = ||\lambda'|| \not= 0$. Let without restriction of generality
$\lambda_0 = 1$. Now choose an orthogonal $(m-1)\times(m-1)$ matrix $U$ such that its first row is given by
${\lambda'}^T$. Define the automorphism $A' \in \Aut(L_m)$ by $A' = \diag(1,U)$
and let $A = {\cal I}_mA'{\cal I}_m^{-1}$ be the corresponding automorphism of $K_m$. Denote the
components of ${\bf B}' = (A \otimes \id_{{\cal S}(n)})({\bf B})$ by $B'_0,\dots,B'_{m-1}$.
We now compute $B'_0$ and $B'_1$. We have $B'_0 = B_0, B'_1 = \sum_{l=1}^{m-1} \lambda_l B_l = -\lambda_0
B_0$. Hence $B'_0+B'_1 = 0$. Note that this sum is the upper left $n \times n$ block of the matrix ${\bf
B}'$. But $A \otimes \id_{{\cal S}(n)} \in {\cal G}_{m,n}$, hence ${\bf B}' \in \Gamma_{m,n}$ by Proposition
\ref{symmetry} and in particular ${\bf B}'$ is PSD. It follows that the first $n$ rows and columns of ${\bf
B}'$ are zero. In particular, we get $B'_{m-1} = 0$, and ${\bf B}' \in E_{m-1} \otimes {\cal S}(n)$.

ii) $|\lambda_0| < ||\lambda'||$. Let without restriction of generality
$||\lambda'||^2 = \lambda_0^2 + 1$. Then there exists a number $\xi$ and a unit length
vector $v \in {\mathbb R}^{m-1}$ such that $\lambda_0 = \sinh\xi, {\lambda'}^T =
\cosh\xi\,v^T$. Let now $A'_1\in \Aut(L_m)$ be a hyperbolic rotation in the $(e_0,e_{m-1})$-plane by the angle
$\xi$, and let $A'_2\in \Aut(L_m)$ be a rotation in the linear subspace spanned by $\{e_1,\dots,e_{m-1}\}$, given by an orthogonal $(m-1) \times (m-1)$ matrix
which has $v$ as its last row.
The last row of the matrix $A' = A'_1A'_2 \in \Aut(L_m)$ is then given by
$(\sinh\xi, \cosh\xi\,v_1,\dots,\cosh\xi\,v_{m-1}) = \lambda^T$. 
Define now the automorphism $A = {\cal I}_mA'{\cal I}_m^{-1} \in \Aut(K_m)$ and let ${\bf B}' = (A \otimes \id_{{\cal S}(n)})({\bf B})$.
The last component of ${\bf B}'$ is then given by $B'_{m-1} = \sum_{l=0}^{m-1} \lambda_l B_l = 0$. Hence ${\bf B}' \in E_{m-1} \otimes {\cal S}(n)$.
\end{proof}

We can hence reduce decomposable elements to elements in spaces with
smaller $n$, and
reducible elements to elements in spaces with smaller $m$. By
using Proposition \ref{symmetry} and applying the lines of reasoning in the proof of Proposition
\ref{reduceimpl} we arrive at
the following proposition.

{\prop A decomposable element of $\Gamma_{m,n}$ is in
$\Sigma_{m,n}$ if and only if its block components are in $\Sigma_{m,n_1}$ and
$\Sigma_{m,n_2}$, respectively. A reducible element of
$\Gamma_{m,n}$ is in $\Sigma_{m,n}$ if and only if its reduction is
in $\Sigma_{m-1,n}$. \qed }

{\corollary \label{reduction} Suppose that $\Gamma_{m,n} =
\Sigma_{m,n}$ for some $n,m$. Then all decomposable elements of
$\Gamma_{m,n+1}$ are in $\Sigma_{m,n+1}$ and all reducible elements
of $\Gamma_{m+1,n}$ are in $\Sigma_{m+1,n}$. \qed }

\medskip

Let us now consider the cases $m = 3$ and $n = 2$.

{\theorem \label{th1} $\Gamma_{3,n} = \Sigma_{3,n}$ for any $n \geq 1$. }

\begin{proof}
The cone $\Gamma_{3,n}$ is the cone of real symmetric PSD
block-Hankel matrices of size $2n \times 2n$. However, such matrices are
known to be separable. This follows from the spectral factorization theorem for quadratic matrix-valued polynomials in one variable \cite{Yakubovitch70}.
\end{proof}

By a similar reasoning applied to complex hermitian block-Hankel matrices we obtain the following result.

{\theorem \label{thH2S} Let $n \geq 1$. A matrix in ${\cal S}(2) \otimes {\cal H}(n)$ is PSD if and only if it is $S_+(2) \otimes H_+(n)$-separable. \qed }

{\theorem \label{th2} $\Gamma_{m,2} = \Sigma_{m,2}$ for $m = 3,\dots,6$. }

\begin{proof}
The space ${\cal S}(2)$ is 3-dimensional. Hence any
quadruple of real symmetric $2\times 2$ matrices is linearly
dependent, and for $m \geq 4$ the space $E_m \otimes {\cal
S}(2)$ consists of reducible elements only. Since $\Gamma_{3,2} =
\Sigma_{3,2}$ by the previous theorem, we have $\Gamma_{m,2} =
\Sigma_{m,2}$ for arbitrary $m = 3,\dots,6$ by repeated application of
Corollary \ref{reduction}.
\end{proof}

{\lemma \label{orbit} Let ${\bf B} \in \Gamma_{4,3}$ be partitioned in symmetric $3 \times 3$ blocks as follows
\begin{equation} \label{Bpartition}
{\bf B} = \left( \begin{array}{cc} B_{11} & B_{12} \\ \overline{B_{12}} & B_{22} \end{array} \right).
\end{equation}
Then either ${\bf B}$ is reducible or decomposable, and hence in $\Sigma_{4,3}$; or in
the orbit of ${\bf B}$ with respect to the action of the group ${\cal G}_{4,3}$ there exists an element 
${\bf B'}$ such that its upper left $3 \times 3$ subblock $B_{11}'$ equals $I_3$. }

\begin{proof} Let the assumptions of the lemma hold. Then we have
${\bf B} \succeq 0$. Hence if a real vector $v \in {\mathbb R}^3$ is in the kernel of one of the
matrices $B_{11},B_{22}$, then it is also in the kernels of $B_{12},\overline{B_{12}}$.

Suppose that neither $B_{11}$ nor $B_{22}$ are of full
rank. 

If the kernels of $B_{11},B_{22}$ have a nontrivial
intersection, then ${\bf B}$ is decomposable, because this intersection
will be in the kernel of all four components of ${\bf B}$.

If the intersection of the kernels is trivial,
then the kernel of $B_{12}$ contains two linearly independent
real vectors. Then the real and imaginary parts of $B_{12}$ must be linearly dependent, because
they are symmetric, sharing a 2-dimensional kernel and of rank not exceeding 1. Hence ${\bf B}$ is reducible.

By Corollary \ref{reduction} and Theorems \ref{th1} and \ref{th2} any
reducible or decomposable element of $\Gamma_{4,3}$ is in $\Sigma_{4,3}$.

Let us now suppose that at least one of the matrices $B_{11},B_{22}$ has full rank and is hence positive definite (PD).
We can assume without loss of generality that this is $B_{11}$, otherwise we pass to a matrix in the same orbit by first applying the automorphism 
$H_3^{\gamma}(\sigma_2) \otimes \id_{{\cal S}(3)} \in {\cal G}_{4,3}$, where $\sigma_2$ is the non-trivial $2 \times 2$ permutation matrix. Then the matrix ${\bf B}' =
(\id_{E_4} \otimes H^{{\mathbb R}}_n((B_{11})^{-1/2}))({\bf B})$ is in the orbit of ${\bf B}$ and we have $B_{11}' = I_3$.
\end{proof}

{\theorem \label{C43} $\Gamma_{4,3} = \Sigma_{4,3}$. }

\begin{proof}
We shall show that any extreme ray of $\Gamma_{4,3}$ is in $\Sigma_{4,3}$.

Let ${\bf B}$ generate an extreme ray of $\Gamma_{4,3}$.
Note that $\Gamma_{4,n}$ is isomorphic to the intersection of the $4n^2$-dimensional cone $H_+(2n)$ with the
$2n(n+1)$-dimensional subspace of matrices consisting of four symmetric $n\times n$ blocks. 
Then the face of ${\bf B}$ in $H_+(6)$ has at most dimension $36 - 24 + 1 = 13$ by Lemma \ref{dimcount}. 
Hence the PSD matrix ${\bf B} \in H_+(6)$ has at most rank 3.

By Lemma \ref{orbit} we can assume without restriction of generality that ${\bf B}$ is partitioned as in (\ref{Bpartition}) with $B_{11} = I_3$.
Since ${\bf B}$ has rank 3, we can factorize it as
\[ \left( \begin{array}{cc} I_3
& B_{12} \\ \overline{B_{12}} & B_{22} \end{array} \right) = \left(
\begin{array}{c} I_3 \\ W \end{array} \right) \left(
\begin{array}{c} I_3 \\ W \end{array} \right)^*
\]
with $W = \overline{B_{12}} = W^T$. It follows that $W\overline{W} = WW^* = B_{22} = (WW^*)^T = \overline{W}W$. Therefore the
real and imaginary parts of $W$ are symmetric and commute. In particular they can be simultaneously
diagonalized by applying an orthogonal transformation $U$. Applying the automorphism $\id_{E_4} \otimes
H^{{\mathbb R}}_n(U)$ to ${\bf B}$, we thus simultaneously diagonalize all four matrices $B_{11},B_{12},\overline{B_{12}},B_{22}$,
and ${\bf B}$ is decomposable and hence separable.

We have proven that all extreme rays of $\Gamma_{4,3}$ are in
$\Sigma_{4,3}$. But then we have $\Gamma_{4,3} = \Sigma_{4,3}$.
\end{proof}

\medskip

In this section we have investigated the relationship between the cones
$\Sigma_{m,n},\Gamma_{m,n}$ for different dimensions $m,n$.
We have defined two properties of elements in $E_m \otimes {\cal S}(n)$,
namely those of being decomposable and reducible.
Our next goal is to show that $\Gamma_{m,n} \subset
\Sigma_{m,n}$ for $n = 3$, $m = 5,6$ (the converse inclusion being trivial).
We have shown in this section that this inclusion is valid for elements
possessing the above-cited properties, provided
the relation $\Gamma_{m,n} = \Sigma_{m,n}$ holds for the respective
cones of smaller dimension. We have proven
this relation for $n \leq 2$; $m \leq 3$; and $m = 4$, $n = 3$. This
allows us to
concentrate on non-decomposable and non-reducible elements in the proofs
of the main results in the next two sections,
which essentially amounts to imposing certain non-degeneracy conditions.

\section{$\Gamma_{5,3} = \Sigma_{5,3}$}

The structure of the proof of this equality resembles that of the
proof of Theorem \ref{C43}. We show that every element of $\Gamma_{5,3}$
generating an extreme ray
is in $\Sigma_{5,3}$.

\medskip

First we derive some properties of real $3 \times 3$ matrices. Let
${\cal A}(n)$ be the space of real skew-symmetric matrices of size $n
\times n$.
The space ${\cal A}(3)$ is isomorphic to ${\mathbb R}^3$. We define an
isomorphism ${\cal V}: {\cal A}(3) \to {\mathbb R}^3$ by
\[ {\cal V}: \left( \begin{array}{ccc} 0 & a_{12} & -a_{31} \\ -a_{12} &
0 & a_{23} \\ a_{31} &
-a_{23} & 0 \end{array} \right) \mapsto \left( \begin{array}{c} a_{23}
\\ a_{31} \\ a_{12}
\end{array} \right).
\]
If $v,w \in {\mathbb R}^3$ are column vectors, then this isomorphism
maps the
skew-symmetric matrix $vw^T - wv^T$ to the cross-product $v \times w$.

Define now a group homomorphism $H^Q: GL_3({\mathbb R}) \to
GL_3({\mathbb R})$ by $[H^Q(C)](v) = {\cal V}(C[{\cal V}^{-1}(v)]C^T)$
for any $v \in {\mathbb R}^3$, $C \in GL_3({\mathbb R})$.
Note that $H^Q(-C) = H^Q(C)$ for all $C \in GL_3({\mathbb R})$. One also easily checks that the induced Lie algebra
homomorphism has a trivial kernel.
Therefore the image of $H^Q$ is the connected component of the identity
matrix $I_3$ and consists of the real $3 \times 3$ matrices
with positive determinant. Direct calculus shows that $H^Q(C) = (\det C)
C^{-T}$.

For 3 column vectors $u,v,w \in {\mathbb R}^3$, let $(u,v,w)$ be the $3
\times 3$ matrix composed of these column vectors.
The following result can be checked by direct computation.

{\lemma \label{detcross} For any three vectors $u,v,w \in {\mathbb R}^3$
we have $(u \times v)
\times (v \times w) = v\,\det(u,v,w)$. Moreover, the following
assertions are equivalent. \\
i) $\det(u,v,w) = 0$, \\
ii) $u,v,w$ are linearly dependent, \\
iii) $u \times v$, $v \times w$, $w \times u$ are linearly dependent, \\
iv) $u \times v$, $v \times w$, $w \times u$ are all proportional. \qed
}

{\corollary \label{signneg} Let $u,v,w \in {\mathbb R}^3$ be linearly
independent vectors. Then $\det(v \times u,w \times v,u \times w) < 0$. }

\begin{proof}
By repeated application of the previous lemma we get
\begin{eqnarray*}
v \times w &=& [(u \times v) \times (v \times w) \,{\det}^{-1}(u,v,w)]
\times [(v \times w) \times (w \times u) \,{\det}^{-1}(v,w,u)] \\ &=&
{\det}^{-2}(u,v,w) [(u \times v) \times (v \times w)] \times [(v \times
w) \times (w \times u)] \\ &=&
{\det}^{-2}(u,v,w) [(v \times w) \det(u \times v,v \times w,w \times u)].
\end{eqnarray*}
It follows that $\det(v \times u,w \times v,u \times w) =
-\det(u \times v,v \times w,w \times u) = -{\det}^2(u,v,w) < 0$.
\end{proof}

We now investigate 2-dimensional traceless linear subspaces ${\cal L}$
in ${\cal S}(3)$. To each such subspace, we will assign a {\sl sign}
$\sigma({\cal L}) \in \{ -1,0,+1 \}$ in the following way.
Let $\{ S_1,S_2,S_3,I_3 \}$ be a basis of the orthogonal complement
${\cal L}^{\perp}$ of ${\cal L}$. Consider the vectors
$v_{12} = {\cal V}(S_1S_2-S_2S_1), v_{23} = {\cal V}(S_2S_3-S_3S_2),
v_{31} = {\cal V}(S_3S_1-S_1S_3) \in {\mathbb R}^3$. Now define the sign
of ${\cal L}$ as
$\sigma({\cal L}) = \sign \det(v_{23},v_{31},v_{12})$.

{\lemma The sign $\sigma({\cal L})$ is well-defined, i.e.\ does not
depend on the choice of the basis of ${\cal L}^{\perp}$. }

\begin{proof}
Let us first remark that ${\cal L}^{\perp}$ has dimension 4 and contains
the identity matrix $I_3$, hence the claimed choice of its basis is
possible.

Let now $\{ S_1,S_2,S_3,I_3 \}$ and $\{ S'_1,S'_2,S'_3,I_3 \}$ be two
bases of ${\cal L}^{\perp}$. Observe that $[A+\alpha I_n,B+\beta I_n] =
[A,B]$ for any two $n \times n$ matrices and any real scalars
$\alpha,\beta$.
Since $\sigma({\cal L})$ depends only on the pairwise commutators of the
basis elements $S_l,S_l'$, we can assume without loss of generality that
these matrices are traceless. Note that then both
$S_l$ and $S_l'$ span the same 3-dimensional space, namely the
orthogonal complement of ${\cal L}$ in the subspace of traceless
symmetric matrices. We hence find a regular $3 \times 3$ matrix $C$
such that $S_l' = \sum_{\alpha = 1}^3 C_{l \alpha} S_{\alpha}$ for $l =
1,2,3$. Here the indexation of $C$ denotes its elements. Define
$v_{\alpha\beta} = {\cal V}([S_{\alpha},S_{\beta}])$, $v_{\alpha\beta}' =
{\cal V}([S_{\alpha}',S_{\beta}'])$, $\alpha,\beta = 1,2,3$. Then we
have by the bilinearity of the matrix commutator that $v_{\alpha\beta}' = \sum_{l,m=1}^n
C_{\alpha l}v_{lm}C_{m\beta}^T$ for $\alpha,\beta = 1,2,3$. It follows that
$(v_{23}',v_{31}',v_{12}') = (v_{23},v_{31},v_{12})(H^Q(C))^T = (\det C)\,(v_{23},v_{31},v_{12})C^{-1}$ and $\det(v_{23}',v_{31}',v_{12}') =
(\det C)^2 \det(v_{23},v_{31},v_{12})$. Thus the signs of
$\det(v_{23}',v_{31}',v_{12}')$ and $\det(v_{23},v_{31},v_{12})$ are equal.
\end{proof}

{\definition We call two 2-dimensional traceless subspaces ${\cal
L},{\cal L}' \subset {\cal S}(3)$ {\sl equivalent} if there exists an
orthogonal $3 \times 3$ matrix $U$ such that
${\cal L}' = \{ USU^T \,|\, S \in {\cal L} \}$. }

{\lemma \label{signeq} For equivalent subspaces ${\cal L},{\cal L}'$ we
have $\sigma({\cal L}) = \sigma({\cal L}')$. }

\begin{proof}
Let $U$ be the orthogonal matrix realizing the equivalence.
Let $\{ S_1,S_2,S_3,I_3 \}$ be a basis of ${\cal L}^{\perp}$ and define
$S_l' = US_lU^T$, $l = 1,2,3$. Then $\{ S'_1,S'_2,S'_3,I_3 \}$ is a
basis of ${{\cal L}'}^{\perp}$. Let further
$v_{\alpha\beta} = {\cal V}([S_{\alpha},S_{\beta}])$, $v_{\alpha\beta}'
= {\cal V}([S_{\alpha}',S_{\beta}'])$, $\alpha,\beta = 1,2,3$. Since we
have $[S_{\alpha}',S_{\beta}'] = U[S_{\alpha},S_{\beta}]U^T$ for all
$\alpha,\beta = 1,2,3$, it follows that $v_{\alpha\beta}' =
[H^Q(U)]v_{\alpha\beta}$ and $(v_{23}',v_{31}',v_{12}') =
[H^Q(U)](v_{23},v_{31},v_{12}) = (\det U)\, U(v_{23},v_{31},v_{12})$. This
finally yields
$\sigma({\cal L}) = \sigma({\cal L}')$.
\end{proof}

Let us now consider symmetric tensors $S_{\alpha\beta\gamma}$ of order 3
in ${\mathbb R}^3$. Since there are 10 independent components, these
tensors form a 10-dimensional real vector space.

{\definition We shall say that a symmetric tensor $S_{\alpha\beta\gamma}$ of
order 3 satisfies the {\sl $\delta$-condition} if
$\sum_{\kappa = 1}^3 (S_{\alpha\beta\kappa}S_{\gamma\eta\kappa} -
S_{\gamma\beta\kappa}S_{\alpha\eta\kappa}) =
\delta_{\gamma\beta}\delta_{\alpha\eta} -
\delta_{\alpha\beta}\delta_{\gamma\eta}$ for
all $\alpha,\beta,\gamma,\eta = 1,2,3$, where $\delta$ is the Kronecker
symbol ($\delta_{\alpha\beta} = 1$ if $\alpha = \beta$ and
$\delta_{\alpha\beta} = 0$ otherwise). }

{\definition Let $S_{\alpha\beta\gamma}$ be a symmetric tensor. Let the
{\sl matrix components} of $S_{\alpha\beta\gamma}$ be three matrices
$S^1,S^2,S^3 \in {\cal S}(3)$ defined elementwise by
$S^l_{\alpha\beta} = S_{\alpha\beta l}$, $\alpha,\beta,l = 1,2,3$. }

{\remark  The $\delta$-condition is equivalent to the condition
$\sum_{l=1}^3 (\det S^l) (S^l)^{-1} = -I_3$, where $S^l$ are the matrix
components of the tensor.
Here the function $A \mapsto (\det A) A^{-1}$ is understood to be
extended by continuity to singular matrices $A \in {\cal S}(3)$. }

{\lemma \label{rotinv} The $\delta$-condition is rotationally invariant,
i.e.\ if $U$ is an orthogonal $3 \times 3$ matrix, then the tensor
$S'_{\alpha\beta\gamma} = \sum_{\eta,\phi,\xi = 1}^3 U_{\alpha\eta}
U_{\beta\phi} U_{\gamma\xi} S_{\eta\phi\xi}$ satisfies the
$\delta$-condition if and only if the tensor $S_{\alpha\beta\gamma}$
does so. }

The proof is by direct calculation using the relation $\sum_{\gamma =
1}^3 U_{\alpha\gamma} U_{\beta\gamma} = \delta_{\alpha\beta}$.

{\lemma Let $S_{\alpha\beta\gamma}$ be a symmetric tensor satisfying the
$\delta$-condition with matrix components $S^l$, $l = 1,2,3$. Then the
matrices $\{ S^1,S^2,S^3,I_3 \}$ are linearly independent. }

\begin{proof}
We proof the lemma from the contrary.
Assume the conditions of the lemma and suppose that there exist
$c_0,c_1,c_2,c_3 \in {\mathbb R}$, not all equal zero,
such that $c_0 I_3 + \sum_{l=1}^3 c_l S^l = 0$. By Lemma \ref{rotinv} we
can assume without loss of generality that
$c_2 = c_3 = 0$. Then there exists $c \not= 0$ such that $S^1 = c I_3$,
or $S_{1\beta\gamma} = c\delta_{\beta\gamma}$ for all $\beta,\gamma =
1,2,3$.
Since $S_{\alpha\beta\gamma}$ satisfies the $\delta$-condition, we have
in particular
$\sum_{\kappa = 1}^3 (S_{11\kappa}S_{22\kappa} - S_{21\kappa}S_{12\kappa})
= \delta_{21}\delta_{12} - \delta_{11}\delta_{22}$. But the left-hand
side of this equation simplifies to 0, whereas the right-hand side
simplifies to $-1$, which leads to a contradiction.
\end{proof}

We now come to a result linking the sign of a traceless 2-dimensional
subspace ${\cal L} \subset {\cal S}(3)$ to symmetric tensors satisfying
the $\delta$-condition.

{\theorem \label{tensorsym} 2-dimensional traceless subspaces of ${\cal
S}(3)$ with sign $-1$ are in correspondence with symmetric tensors
satisfying the $\delta$-condition. Namely, if ${\cal L}$ is such a subspace,
then there exists a symmetric tensor $S_{\alpha\beta\gamma}$ with matrix components $S^l$, $l = 1,2,3$ satisfying
the $\delta$-condition such that $\{ S^1,S^2,S^3,I_3 \}$ is a basis of
the space ${\cal L}^{\perp}$. Conversely, if
$S_{\alpha\beta\gamma}$ is such a tensor, then the orthogonal complement of the linear span of the
set $\{ S^1,S^2,S^3,I_3 \}$ has sign $-1$. }

\begin{proof}
Let $e_1,e_2,e_3$ be the canonical orthonormal basis vectors of
${\mathbb R}^3$.

Let $S_{\alpha\beta\gamma}$ be a tensor with matrix components $S^l$,
$l=1,2,3$ satisfying the $\delta$-condition. Then we have $S^1S^2-S^2S^1
= e_2e_1^T-e_1e_2^T$ and
${\cal V}([S^1,S^2]) = e_2 \times e_1 = -e_3$. Similarly, ${\cal
V}([S^2,S^3]) = -e_1$ and ${\cal V}([S^3,S^1]) = -e_2$. Let ${\cal L}$
be the orthogonal complement of the linear span of the set $\{
S^1,S^2,S^3,I_3 \}$.
Then by definition $\sigma({\cal L}) = \sign \det(-e_1,-e_2,-e_3) = -1$,
which proves the second part of the theorem.

Let us prove the first part. Let ${\cal L} \subset {\cal S}(3)$ be a
2-dimensional traceless subspace with sign $\sigma({\cal L}) = -1$.
Any point in ${\cal S}(3)$ can be viewed as a homogeneous quadratic form
on ${\mathbb R}^3$, or equivalently, as a quadratic map from ${\mathbb
RP}^2$ to ${\mathbb R}$.
We are interested in the number of points in ${\mathbb RP}^2$ which are
mapped to zero by all elements of a  ${\cal L}$. Denote $N_{\cal L} = \{
x \in {\mathbb R}^3 \,|\, x^TSx = 0 \ \forall\ S \in {\cal L} \}$.

The determinant as a scalar function on ${\cal L}$ is odd and hence
possesses zeros on ${\cal L} \setminus \{ 0 \}$. Since the matrices in
${\cal L}$ are traceless, we can find a matrix in ${\cal L}$ with eigenvalues $-1,0,+1$.
By conjugation with an appropriate orthogonal matrix $U$ we can
transform it to the matrix
\[ N_1 = \left( \begin{array}{ccc} 0 & 1 & 0 \\ 1 & 0 & 0 \\ 0 & 0 & 0
\end{array} \right).
\]
Let ${\cal L}' = \{ USU^T \,|\, S \in {\cal L} \}$. Then the subspaces
${\cal L}$ and ${\cal L}'$ are equivalent and we have $N_1 \in {\cal L}'$.
A generator the orthogonal complement of $N_1$ in ${\cal L}'$ will be of
the form
\[ N_2 = \left( \begin{array}{ccc} a & 0 & c \\ 0 & b & d \\ c & d &
-a-b \end{array} \right)
\]
for some numbers $a,b,c,d \in {\mathbb R}$, which do not equal zero
simultaneously. Note that by Lemma \ref{signeq} we have $\sigma({\cal
L}') = -1$.
Let us first treat several degenerate cases.

\smallskip

1. $a = b = 0$. In this case a basis of ${{\cal L}'}^{\perp}$ is given
by $\{ S_1,S_2,S_3,I_3 \}$ with
\[ S_1 = \left( \begin{array}{ccc} 1 & 0 & 0 \\ 0 & 0 & 0 \\ 0 & 0 & 0
\end{array} \right),\ S_2 = \left( \begin{array}{ccc} 0 & 0 & 0 \\ 0 & 1
& 0 \\ 0 & 0 & 0 \end{array} \right),\
S_3 = \left( \begin{array}{ccc} 0 & 0 & d \\ 0 & 0 & -c \\ d & -c & 0
\end{array} \right).
\]
We get $\sigma({\cal L}') = 0$, because $[S_1,S_2] = 0$.

\smallskip

2. $a = -b \not= 0$. In this case a basis of ${{\cal L}'}^{\perp}$ is
given by $\{ S_1,S_2,S_3,I_3 \}$ with
\[ S_1 = \left( \begin{array}{ccc} 0 & 0 & 0 \\ 0 & 0 & 0 \\ 0 & 0 & 1
\end{array} \right),\ S_2 = \left( \begin{array}{ccc} -ad & 0 & -cd \\ 0
& ad & a^2+c^2 \\ -cd & a^2+c^2 & 0 \end{array} \right),\
S_3 = \left( \begin{array}{ccc} -ac & 0 & a^2+d^2 \\ 0 & ac & -cd \\
a^2+d^2 & -cd & 0 \end{array} \right).
\]
We get $\det(v_{23},v_{31},v_{12}) = a^4(a^2+c^2+d^2)^2 > 0$ and
$\sigma({\cal L}') = 1$. Here $v_{\alpha\beta} = {\cal
V}([S_{\alpha},S_{\beta}])$.

\smallskip

Hence these two cases do not satisfy the conditions of the theorem and
we can assume $a+b \not= 0$.

The set $N_{{\cal L}'}$ is given by those vectors $x = (x_1,x_2,x_3)^T
\in {\mathbb R}^3$ that satisfy $x^TN_1x = x^TN_2x = 0$. In particular,
$x$ must satisfy $x_1x_2 = 0$. If we define
\[ N_{21} = \left( \begin{array}{cc} a & c \\ c & -a-b \end{array}
\right), \quad N_{22} = \left( \begin{array}{cc} b & d \\ d & -a-b
\end{array} \right),
\]
then we get
\begin{equation} \label{N_L}
N_{{\cal L}'} = \left\{ \left( \begin{array}{c} x_1 \\ 0 \\ x_3
\end{array} \right) \,|\, (x_1\ x_3) N_{21} \left( \begin{array}{c} x_1
\\ x_3 \end{array} \right) = 0 \right\} \cup
\left\{ \left( \begin{array}{c} 0 \\ x_2 \\ x_3 \end{array} \right)
\,|\, (x_2\ x_3) N_{22} \left( \begin{array}{c} x_2 \\ x_3 \end{array}
\right) = 0 \right\}.
\end{equation}
Note that these two sets have a trivial intersection due to the
condition $a+b \not= 0$.
Let $D_1 = \det N_{21}$, $D_2 = \det N_{22}$. We have $D_1+D_2 =
-(a+b)^2-c^2-d^2 < 0$ and at least one of these two determinants is
negative. Let this be $D_1$ without restriction of generality (otherwise
we exchange $a,b$; $c,d$; $x_1,x_2$; and $N_{21},N_{22}$).

A basis of ${{\cal L}'}^{\perp}$ is then given by $\{ S_1,S_2,S_3,I_3
\}$ with
\begin{eqnarray*}
S_1 &=& \left( \begin{array}{ccc} (a+b)^2 & 0 & (a+b)(c+\sqrt{-D_1}) \\
0 & 0 & 0 \\ (a+b)(c+\sqrt{-D_1}) & 0 & (c+\sqrt{-D_1})^2 \end{array}
\right),\\
S_2 &=& \left( \begin{array}{ccc} (a+b)^2 & 0 & (a+b)(c-\sqrt{-D_1}) \\
0 & 0 & 0 \\ (a+b)(c-\sqrt{-D_1}) & 0 & (c-\sqrt{-D_1})^2 \end{array}
\right),\\
S_3 &=& \left( \begin{array}{ccc} -ad & 0 & -cd \\ 0 & -bd &
a^2+ab+b^2+c^2 \\ -cd & a^2+ab+b^2+c^2 & (a+b)d \end{array} \right).
\end{eqnarray*}
We then get $\det(v_{23},v_{31},v_{12}) =
-4ab(-D_1)(a+b)^4(a^2+ab+b^2+c^2)^2$ and $\sigma({\cal L}') = -\sign
ab$. Since $\sigma({\cal L}') = -1$, we get $ab > 0$. Since we can multiply $N_2$ by $-1$, we
can assume without restriction of generality that $a,b > 0$. We have
$D_2 = -ab-b^2-d^2 < 0$. Hence both $N_{21}$ and $N_{22}$ have negative
determinant.
But the number of linearly independent vectors in the two subsets on the
right-hand side of (\ref{N_L}) depends precisely on the sign of these
determinants. Namely, if $N \not = 0$
is a real symmetric $2 \times 2$ matrix with negative determinant, then
the equation $y^TNy = 0$ has two linearly independent solutions $y \in
{\mathbb R}^2$.
Thus $N_{{\cal L}'}$ contains 4 linear 1-dimensional subspaces of
${\mathbb R}^3$.

Let $x^l$, $l = 1,\dots,4$ be generators of these subspaces, defined by
the columns of the matrix
\[ X = \left( \begin{array}{cccc}
\sqrt{\frac{a}{b}}\frac{a+b}{\sqrt{-D_1}+c} &
-\sqrt{\frac{a}{b}}\frac{a+b}{\sqrt{-D_1}-c} & 0 & 0 \\
0 & 0 & -\sqrt{\frac{b}{a}}\frac{a+b}{\sqrt{-D_2}+d} &
\sqrt{\frac{b}{a}}\frac{a+b}{\sqrt{-D_2}-d} \\
\sqrt{\frac{a}{b}} & \sqrt{\frac{a}{b}} & -\sqrt{\frac{b}{a}} &
-\sqrt{\frac{b}{a}} \end{array} \right).
\]
Each 3 of these vectors span the whole space ${\mathbb R}^3$ and
\begin{equation} \label{xsc}
\langle x^{\alpha},x^{\beta} \rangle = -1\qquad \forall\ \alpha,\beta =
1,2,3,4,\ \alpha \not= \beta.
\end{equation}
The rank 1 matrices $x^l(x^l)^T \in {\cal S}(3)$ are linearly independent and hence form
a basis of the space ${{\cal L}'}^{\perp}$.

Let us define symmetric matrices $X^l$, $l = 1,\dots,4$ by $X^l =
x^l(x^l)^T - I_3$, and positive numbers $c_l = \frac{1}{1 + |x^l|^2}$.
Consider the symmetric $4 \times 4$ matrix $X^TX$. It has rank 3 and is hence singular.
Its diagonal elements are given by $|x^l|^2$, while its off-diagonal elements equal $-1$ by (\ref{xsc}).
Setting the determinant to zero and transforming the obtained relation, we get
$\sum_{l=1}^4 c_l = 1$.
Then using (\ref{xsc}) one checks that
$\langle \sum_{l = 1}^4 c_l x^l, x^p \rangle = \langle \sum_{l=1}^4 c_l X^l, x^p(x^p)^T \rangle = 0$ for all $p = 1,\dots,4$.
This implies $\sum_{l = 1}^4 c_l x^l = \sum_{l=1}^4 c_l X^l = 0$. Hence the matrices
$X^l$ span a 3-dimensional subspace ${\cal L}'_3 \subset {{\cal
L}'}^{\perp}$, which does not contain the matrix $I_3$,
and there exists a unique linear mapping ${\cal X}: {\mathbb R}^3 \to
{\cal L}'_3$ such that ${\cal X}(x^l) = X^l$ for all $l = 1,\dots,4$.
Moreover, we have again by (\ref{xsc}) that
$X^{\alpha}x^{\beta} = X^{\beta}x^{\alpha}$ and $X^{\alpha}X^{\beta} -
X^{\beta}X^{\alpha} = x^{\beta}(x^{\alpha})^T - x^{\alpha}(x^{\beta})^T$
for all $\alpha,\beta = 1,2,3,4$.

By bilinearity these relations extend to ${\cal X}(v)w = {\cal X}(w)v$
and ${\cal X}(v){\cal X}(w) - {\cal X}(w){\cal X}(v) = wv^T - vw^T$ for
all $v,w \in {\mathbb R}^3$. Let now
$S^l = {\cal X}(e_l)$, $l = 1,2,3$. Then by the first relation the
matrices $S^l$ are the matrix components of some symmetric tensor
$S'_{\alpha\beta\gamma}$, and by the second relation this tensor
satisfies the $\delta$-condition. Thus $S'_{\alpha\beta\gamma}$
satisfies the assertions of the theorem with respect to the subspace
${\cal L}'$. For reference, the components of this tensor equal
$S'_{111} = -\frac{2c}{\sqrt{ab}}$, $S'_{222} = \frac{2d}{\sqrt{ab}}$,
$S'_{113} = \sqrt{\frac{a}{b}}$, $S'_{223} = -\sqrt{\frac{b}{a}}$,
$S'_{333} = \sqrt{\frac{a}{b}} - \sqrt{\frac{b}{a}}$, all other
independent components being zero.

Finally, define the tensor $S_{\alpha\beta\gamma} = \sum_{\eta,\phi,\xi
= 1}^3 U_{\alpha\eta}^T U_{\beta\phi}^T U_{\gamma\xi}^T
S'_{\eta\phi\xi}$, where $U$ is the orthogonal matrix realizing the
equivalence
of the subspaces ${\cal L}$ and ${\cal L}'$. By Lemma \ref{rotinv} this
tensor satisfies the assertions of the theorem with respect to ${\cal L}$.
\end{proof}


Finally we provide the following auxiliary result about quaternion vectors.

{\lemma \label{quat3} Let $v \in {\mathbb H}^3$ and $p \in {\mathbb H}$,
$p \not= 0$ be given. Then there exists
a unit quaternion $q$ such that $Re(vq\bar p) = 0$. Similarly, there exists
a unit quaternion $q'$ such that $Re(q'v\bar p) = 0$. }

\begin{proof}
The equation $Re(vq\bar p) = 0$ (or $Re(q'v\bar p) = 0$) amounts to 3 linear
relations on the 4 real components of $q$ (or $q'$). Hence there exists a
non-zero solution,
which can be normalized.
\end{proof}

In particular, we can multiply any vector in ${\mathbb H}^3$ by a unit
quaternion from the
left or from the right such that a given component (e.g.\ the real part)
vanishes.

\medskip

Let us return to the cones $\Gamma_{5,n}$ and $\Sigma_{5,n}$.
Our goal is to show that these cones are equal for $n = 3$.

{\lemma \label{orbit5} Any matrix ${\bf B} \in \Gamma_{5,3}$ is either
reducible or decomposable, and hence in $\Sigma_{5,3}$, or
in the orbit of ${\bf B}$ with respect to the action of the group ${\cal G}_{5,3}$
there exists an element ${\bf B'}$ such that its upper left $3 \times 3$ block $B_{11}'$ equals
$I_3$. \qed }

The proof is similar to the proof of Lemma \ref{orbit}, with a
reference to Theorem \ref{C43} instead of Theorem \ref{th1}.

\medskip

Let now ${\bf B}$ be the generator of an extreme ray of $\Gamma_{5,3}$. Our goal is to show that ${\bf B} \in
\Sigma_{5,3}$. Let us apply the dimensional argument Lemma \ref{dimcount}.

{\corollary \label{rank45} Let ${\bf B}$ generate 
an extreme ray of $\Gamma_{5,3}$. Then ${\bf B}$ has at most
rank 4. }

\begin{proof}
Any element ${\bf B} \in \Gamma_{5,n}$ satisfies ${\bf B} \succeq 0$.
The space ${\cal Q}(6)$ has dimension 66, whereas its subspace $E_5
\otimes {\cal S}(n)$
has dimension 30. Hence the face in $Q_+(6)$ of the extreme ray ${\bf
B}$ has at most dimension 37 by Lemma \ref{dimcount}. Since the space ${\cal Q}(5)$ has
dimension 45, a rank 5 element cannot be extremal (cf.\ Proposition \ref{quatface} in the appendix).
\end{proof}

{\lemma \label{14} Let ${\bf B}$ generate an
extreme ray of $\Gamma_{5,3}$. Then ${\bf B} \in \Sigma_{5,3}$. }

\begin{proof}
Assume the conditions of the lemma.
Let ${\bf B}$ be partitioned in 4 symmetric $n \times n$ blocks as in (\ref{Bpartition}).
By Lemma \ref{orbit5} and Proposition \ref{symmetry} we can assume
without loss of generality that $B_{11} = I_3$ and ${\bf B}$ is not reducible.

Denote the $3\times3$ matrix $\overline{B_{12}}$ by $W$.
Then $B_{22} - WW^*$ is a PSD matrix, which by Corollary \ref{rank45}
has rank at most 1 (cf.\ Lemma \ref{quat_schur_compl} in the appendix). 
Hence there exists a quaternionic vector $Z \in
{\mathbb H}^3$
such that $ZZ^* + WW^* = B_{22}$. Here $Z$ can also be the zero vector.
The vector $Z$ is determined up to multiplication with a scalar of unit
norm from the right.
By Lemma \ref{quat3}, we can find a unit quaternion $q$ such
that $Zq$ has a zero $k$-component. Let us hence assume without loss of
generality
that $Z$ has a zero $k$-component.

Let us factorize ${\bf B}$. We have
\[ {\bf B} = \left( \begin{array}{cc} I_3 & W^* \\ W & B_{22}
\end{array} \right) =
\left( \begin{array}{cc} I & 0 \\ W & Z \end{array} \right)
\left( \begin{array}{cc} I & 0 \\ W & Z \end{array} \right)^*
\]
Note that $W$ is symmetric and with zero $k$-component. We have as well
that $WW^* + ZZ^* = B_{22}$ is symmetric and real. This implies
\begin{eqnarray} \label{lemcond}
W_iW_r - W_rW_i + Z_iZ_r^T - Z_rZ_i^T &=& 0, \nonumber\\
W_jW_r - W_rW_j + Z_jZ_r^T - Z_rZ_j^T &=& 0, \nonumber\\
W_jW_i - W_iW_j + Z_jZ_i^T - Z_iZ_j^T &=& 0,
\end{eqnarray}
where the indices denote the corresponding components of the matrix $W$
and the vector $Z$.
Let us denote the row vector $Z^T$ by $Y$. Then we have $Z^* = \overline Y$.
Since ${\bf B}$ is not reducible, the four matrices $W_r,W_i,W_j,I_3$ are linearly independent.
Denote the orthogonal complement of their linear span by ${\cal L}$.
We shall consider several cases.

\medskip

1. The components of $Z$ are linearly independent.

Condition (\ref{lemcond}) and Corollary \ref{signneg} imply that ${\cal
L}$ has negative sign.
By Theorem \ref{tensorsym}, there exists a symmetric third order tensor
$S_{\alpha\beta\gamma}$ fulfilling the $\delta$-condition whose matrix
components span a 3-dimensional subspace
${\cal L}_3 \subset {\cal L}^{\perp}$. This subspace is transversal to
$I_3$, therefore there exist real numbers $x_r,x_i,x_j$ such that the
matrices $W_r-x_r I_3,W_i-x_i I_3,W_j-x_j I_3$ span
${\cal L}_3$.

Let us now consider the symmetric completion of the $3\times4$
quaternionic matrix $(W\ Z)$
to the square matrix
\[ N = \left( \begin{array}{cc} W & Z \\ Y & x \end{array}
\right),
\]
where $x = x_r + ix_i + jx_j$ is a quaternion with zero $k$-component.
The matrix $N$ has zero $k$-component and mutually commuting real, $i$- and
$j$-component.  Namely, the commutation condition amounts to the equations
\begin{eqnarray*}
[W_{\alpha},W_{\beta}] &=& Z_{\beta}Y_{\alpha} - Z_{\alpha}Y_{\beta}, \\
W_{\alpha}Z_{\beta} + Z_{\alpha}x_{\beta} &=& W_{\beta}Z_{\alpha} +
Z_{\beta}x_{\alpha}, \\
Y_{\alpha}Z_{\beta} + x_{\alpha}x_{\beta} &=& Y_{\beta}Z_{\alpha} +
x_{\beta}x_{\alpha},
\end{eqnarray*}
for every pair of indices $\alpha,\beta$ from $\{r,i,j\}$.
The first set of equations is precisely (\ref{lemcond}), the
last set of equations is satisfied for any $x$
because $Y_{\alpha} = Z_{\alpha}^T$ for all indices $\alpha$.
These second set follows from the symmetricity of the tensor
$S_{\alpha\beta\gamma}$.

Since the matrices $N_r,N_i,N_j$ mutually commute and are symmetric,
they share a
common orthonormal set of real eigenvectors $v_1,v_2,v_3,v_4 \in
{\mathbb R}^4$. This set is also a
set of eigenvectors to the quaternionic matrix $N$, with
eigenvalues $q_1,q_2,q_3,q_4$ that have zero $k$-component. Hence we have
$(W\ Z)v_l = q_l(I_3\ 0)v_l = q_l\tilde
v_l$ for all
$l \in \{1,2,3,4\}$. Here $\tilde v_l \in {\mathbb R}^3$
consists of the first three components of $v_l$.
If we decompose the identity matrix $I_4$ as
$\sum_{l=1}^4 v_lv_l^T$, we can rewrite ${\bf B}$ as
\[ {\bf B} = \sum_{l = 1}^4\ \left( \begin{array}{cc} I_3 & 0 \\ W
& Z \end{array} \right)
v_lv_l^T \left( \begin{array}{cc} I_3 & 0 \\ W & Z
\end{array} \right)^* =
\sum_{l = 1}^4\ \left( \begin{array}{cc} 1 & \overline{q_l} \\
q_l & |q_l|^2
\end{array} \right) \otimes (\tilde v_l \tilde v_l^T ).
\]
This matrix is thus $K_5 \otimes S_+(n)$-separable, and ${\bf B} \in
\Sigma_{5,3}$.

\medskip

2. The components of $Z$ are linearly dependent.

Let $z \in {\mathbb R}^3$ be a nonzero vector in the orthogonal
complement to the subspace spanned by the components of $Z$.
Then by (\ref{lemcond}) we have $[W_r,W_i]z = [W_i,W_j]z = [W_j,W_r]z = 0$.
By bilinearity we then have $[U,V]z = 0$ for any two matrices $U,V \in
{\cal L}^{\perp}$.
Define the linear subspace ${\cal S}_W = \{ U \in {\cal L}^{\perp} \,|\,
Uz = 0 \}$ and let ${\cal S}_z \subset {\mathbb R}^3$ be the linear
span of the vectors $z = I_3z,W_rz,W_iz,W_jz$. Obviously we have
\begin{equation} \label{dimcond}
\dim {\cal S}_W + \dim {\cal S}_z = 4.
\end{equation}
For any matrices $U \in {\cal S}_W$ and $V \in {\cal L}^{\perp}$ we have
$U(Vz) = V(Uz) = 0$. It follows that for any matrix $U \in {\cal S}_W$
and any vector $y \in {\cal S}_z$ we have $Uy = 0$. If now $\dim {\cal
S}_z = 3$, then ${\cal S}_W$ can consist only of the zero matrix, which
contradicts
(\ref{dimcond}). Likewise, if $\dim {\cal S}_z = 2$, then all matrices
in ${\cal S}_W$ are proportional to the rank 1 matrix $vv^T$, where $v$
is the orthogonal complement to ${\cal S}_z$.
Hence $\dim {\cal S}_W \leq 1$, which also contradicts (\ref{dimcond}).
Therefore $\dim {\cal S}_z = 1$, and $z$ is a common eigenvector of all
matrices in ${\cal L}^{\perp}$.
It follows that $z$ is an eigenvector of $W,W^*$, and also $ZZ^*$,
because $Z^*z = 0$. Thus $z$ is a common eigenvector of the 5 matrices
$B_{11}=I_3,W_r,W_i,W_j,B_{22}$ and these matrices are linearly dependent.
But then the components of ${\bf B}$ are also linearly dependent, which
contradicts the non-reducibility of ${\bf B}$.

This completes the proof of Lemma \ref{14}.
\end{proof}

We obtain the following theorem.

{\theorem \label{C53} $\Gamma_{5,3} = \Sigma_{5,3}$. \qed }

\section{$\Gamma_{6,3} = \Sigma_{6,3}$}

In this section we prove the analogue of Theorem \ref{C53} for dimension $m = 6$, following essentially the
same line of reasoning. First we provide some auxiliary results.

{\lemma \label{SU2SO4} Let $v,w \in {\mathbb H}^3$ be quaternionic column vectors of length $3$ such that $Re\,
vv^* = Re\, ww^*$. Then there exist unit norm quaternions $h,h'$ such that $h'v = wh$.}

\begin{proof}
Let the assumptions of the lemma hold. The quaternions ${\mathbb H}$ can be considered as a 4-dimensional
real vector space, equipped with the Euclidean scalar product $\langle a,b \rangle = Re\, a\bar b$. Let ${\cal
W}: {\mathbb H} \to {\mathbb R}^4$ be the corresponding isomorphism. We have for any pair of indices
$\alpha,\beta = 1,2,3$ that $Re\, v_{\alpha} \overline{v_{\beta}} = Re\, w_{\alpha} \overline{w_{\beta}}$. We then get
$\langle {\cal W}(v_{\alpha}),{\cal W}(v_{\beta}) \rangle = \langle {\cal W}(w_{\alpha}),{\cal W}(w_{\beta})
\rangle$. Therefore there exists an orthogonal $4 \times 4$-matrix $U$ such that ${\cal W}(w_{\alpha}) =
U{\cal W}(v_{\alpha})$ for all $\alpha = 1,2,3$. Since $v,w$ have only 3 components, this matrix can be
chosen to have determinant 1. But the group generated by multiplication of ${\mathbb H}$ by unit quaternions
from the left and from the right is known to be isomorphic to the special orthogonal group $SO(4)$. This can
be easily checked by comparing the Lie algebras of these groups. Hence there exist unit norm quaternions
$h,h'$ such that $w = h'vh^{-1}$.
\end{proof}

We provide the following lemma on the characteristic polynomial of a $3\times3$ hermitian quaternionic
matrix.

{\lemma \label{charpoly} Let $A$ be a hermitian quaternionic $3\times3$ matrix with elements $a_{\alpha\beta}
= a_{\alpha\beta}^r + ia_{\alpha\beta}^i + ja_{\alpha\beta}^j + ka_{\alpha\beta}^k$, $\alpha,\beta = 1,2,3$. Then the coefficient
$c_1$ of its characteristic polynomial $\lambda^3 + c_2\lambda^2 + c_1\lambda + c_0$ is given by
\begin{eqnarray*}
c_1 &=& -(a^k_{23})^2+a^r_{11}a^r_{22}-(a^i_{13})^2+a^r_{22}a^r_{33}+a^r_{11}a^r_{33}
-(a^r_{12})^2-(a^j_{13})^2-(a^r_{13})^2-(a^j_{23})^2 \\
&& -(a^j_{12})^2-(a^r_{23})^2-(a^i_{12})^2-(a^k_{12})^2-(a^i_{23})^2-(a^k_{13})^2. \quad \qed
\end{eqnarray*}
}
The proof is by direct calculation (cf.\ Proposition \ref{charpolquat} in the appendix).

\medskip

The partial transposition in the space $E_6 \otimes {\cal S}(n)$ amounts
to quaternionic conjugation.
Hence the cone $\Gamma_{6,m}$ is composed of those ${\bf B} \in E_6 \otimes {\cal S}(3)$ that satisfy the
linear matrix inequalities
\begin{equation} \label{Q6}
{\bf B} = \left( \begin{array}{cc} B_{11} & B_{12} \\
\overline{B_{12 }} & B_{22} \end{array} \right) \succeq 0, \qquad
\overline{{\bf B}} = \left( \begin{array}{cc} B_{11} &
\overline{B_{12 }} \\
B_{12} & B_{22} \end{array} \right) \succeq 0.
\end{equation}
Here ${\bf B}$ is partitioned in 4 symmetric $3 \times 3$ blocks as in (\ref{Bpartition}).

{\lemma \label{orbit6} Let ${\bf B} \in
\Gamma_{6,3}$. Then either
${\bf B}$ is reducible or decomposable, and hence in $\Sigma_{6,3}$, or
in the orbit of
${\bf B}$ with respect to the action of the group ${\cal G}_{6,3}$
there exists an element ${\bf B'}$ whose upper left $3 \times 3$ block $B_{11}'$ equals
$I_3$. \qed }

The proof is similar to the proof of Lemma \ref{orbit}, with a reference
to Theorem \ref{C53} instead
of Theorem \ref{th1}.

\medskip

We now proceed as in the previous section and look for elements ${\bf B} \in \Gamma_{6,3}$ which are not in $\Sigma_{6,3}$.
Let ${\bf B}$ be partitioned as in (\ref{Q6}).
By Lemma \ref{orbit6} and Proposition \ref{symmetry} we can assume without loss of generality that ${\bf B}$ is not reducible and that $B_{11} = I_3$.

We can factorize the PSD matrices ${\bf B},\overline{{\bf B}}$ as
\begin{equation} \label{factorization}
{\bf B} = \left( \begin{array}{cc} I & 0 \\ W & Z \end{array} \right)
\left( \begin{array}{cc} I & 0 \\ W & Z \end{array} \right)^*,
\quad \overline{{\bf B}} = \left( \begin{array}{cc} I & 0 \\ \overline{W} & Y^*
\end{array} \right) \left( \begin{array}{cc} I & 0 \\ \overline{W} & Y^*
\end{array} \right)^*,
\end{equation}
where $W = \overline{B_{12 }}$, and $Y,Z$ are quaternionic matrices of
appropriate size (cf.\ Lemma \ref{quat_schur_compl} in the appendix). 
Note that $W$ and $W\overline{W} + ZZ^* = \overline{W}W + Y^*Y = B_{22}$ are symmetric matrices.

As in the previous section, denote the components of $W,Z,Y$ by adding a
corresponding index. The imaginary part of $W\overline{W} + ZZ^*$ is zero,
which yields
\begin{eqnarray} \label{commrel6}
\,[W_i,W_r] + [W_k,W_j] + Z_iZ_r^T - Z_rZ_i^T + Z_kZ_j^T - Z_jZ_k^T &=&
0, \nonumber\\
\,[W_j,W_r] + [W_i,W_k] + Z_jZ_r^T - Z_rZ_j^T + Z_iZ_k^T - Z_kZ_i^T &=&
0, \nonumber\\
\,[W_k,W_r] + [W_j,W_i] + Z_kZ_r^T - Z_rZ_k^T + Z_jZ_i^T - Z_iZ_j^T &=& 0.
\end{eqnarray}
In a similar way, setting the imaginary part of $\overline{W}W + Y^*Y$ to
zero, we obtain
\begin{eqnarray} \label{commrel62}
-[W_i,W_r] + [W_k,W_j] - Y^T_iY_r + Y^T_rY_i + Y^T_kY_j - Y^T_jY_k &=&
0, \nonumber\\
-[W_j,W_r] + [W_i,W_k] - Y^T_jY_r + Y^T_rY_j + Y^T_iY_k - Y^T_kY_i &=&
0, \nonumber\\
-[W_k,W_r] + [W_j,W_i] - Y^T_kY_r + Y^T_rY_k + Y^T_jY_i - Y^T_iY_j &=& 0.
\end{eqnarray}

Let us apply the dimensional argument Lemma \ref{dimcount}. Referring to the notations of this lemma, we
choose $K$ to be the 132-dimensional direct product $Q_+(6)\times Q_+(6)$. Its faces are isomorphic to direct
products $Q_+(l_1) \times Q_+(l_2)$ with $l_1 \leq 6$, $l_2 \leq 6$ (cf.\ Proposition \ref{quatface} in the appendix). The subspace $L$ consists of all pairs
$(B,\overline B) \in {\cal Q}(6) \times {\cal Q}(6)$ with $B$ consisting of four symmetric $3\times3$ blocks. This
space is 36-dimensional. The intersection $K' = K \cap L$ is isomorphic to the cone of blockwise symmetric
PPT matrices in $Q_+(6)$, i.e.\ to $\Gamma_{6,3}$. By Lemma \ref{dimcount}, an extremal ray of $\Gamma_{6,3}$
lies in a face of $K$ that has at most dimension 97. Now note that the space ${\cal Q}(5) \times {\cal Q}(6)$
has dimension $45+66 = 111$. Hence if $B$ or $\overline B$ are of full rank, then the corresponding conjugate can
be at most of rank 4. We get the following corollary.

{\corollary \label{cor14} Let ${\bf B}$, factorized as in (\ref{factorization}), define an extreme ray of $\Gamma_{6,3}$. 
Then either $\rk\,ZZ^* \leq 1$ or $\rk\,Y^*Y \leq 1$ or $\rk\,ZZ^* = \rk\,Y^*Y = 2$. \qed }

\smallskip

We treat these cases separately. Let the assumptions of Corollary
\ref{cor14} hold and suppose that ${\bf B}$ is not reducible.

\medskip

1. The matrix $ZZ^*$ has at most rank 1, i.e.\ $Z$ is a quaternionic
vector, which can be zero.

\medskip

1.1. The four components of $Z$ lie in a 2-dimensional linear subspace
of ${\mathbb R}^3$.

Then there exists a non-zero vector $z \in {\mathbb R}^3$ such that $Z^*z = 0$. We have $W\overline{W} - \overline{W}W =
W(2W_r - W) - (2W_r - W)W = 2[W,W_r]$. It follows that $||Yz||^2 = z^*Y^*Yz = z^*(2[W,W_r] + ZZ^*)z =
2z^T[W,W_r]z = 0$. The last equality follows from the fact that $z$ is real and hence commutes with
quaternions, and that $[W,W_r]$ is skew-symmetric. Therefore we also have $Yz = 0$. Multiplying
(\ref{commrel6}) and (\ref{commrel62}) by $z$ from the right and combining, we obtain
$[W_{\alpha},W_{\beta}]z = 0$ for all index pairs $\alpha,\beta = r,i,j,k$. By bilinearity we then have
$[U,V]z = 0$ for any two real symmetric matrices $U,V$ in the linear span ${\cal L}$ of
$I_3,W_r,W_i,W_j,W_k$. Note that these 5 matrices are linearly independent, because ${\bf B}$ is not
reducible. Define the linear subspace ${\cal S}_W = \{ U \in {\cal L} \,|\, Uz = 0 \}$ and let ${\cal S}_z
\subset {\mathbb R}^3$ be the linear span of the vectors $z = I_3z,W_rz,W_iz,W_jz,W_kz$. Obviously we have
\begin{equation} \label{dimcond2}
\dim {\cal S}_W + \dim {\cal S}_z = 5.
\end{equation}
For any matrices $U \in {\cal S}_W$ and $V \in {\cal L}$ we have $U(Vz) = V(Uz) = 0$. It follows that for any
matrix $U \in {\cal S}_W$ and any vector $y \in {\cal S}_z$ we have $Uy = 0$. By the same reasoning as in the
previous section we can lead the cases $\dim {\cal S}_z = 3$ and $\dim {\cal S}_z = 2$ to a contradiction
with (\ref{dimcond2}). Therefore $\dim {\cal S}_z = 1$, and $z$ is a common eigenvector of all matrices in
${\cal L}$. This yields two independent linear conditions on the matrices in ${\cal L}$, but ${\cal L}$ has
only codimension 1 in ${\cal S}(3)$. Therefore the considered case is not possible.

\medskip

1.2. The components of $Z$ span the whole space ${\mathbb R}^3$.

We have $Y^*Y = WW^* - W^*W + ZZ^*$. It follows that $Re(Y^*Y) = Re(ZZ^*)$ and $Im(Y^*Y) = -Im(W^*W)$. The
coefficient $c_1$ of the characteristic polynomial of $Y^*Y$ is therefore a homogeneous polynomial of degree
4 in the components of the elements of $W,Z$. By Lemma \ref{charpoly} the parts that depend on $W$ and on $Z$
are separated. It is not hard to verify that $Z$ enters only in the form of squares of the elements of the
imaginary part of $ZZ^*$. If we replace these elements by corresponding terms depending on $W$ using the
relation $Im(ZZ^*) = -Im(WW^*)$, then the polynomial simplifies to zero. Hence the derivative of the
characteristic polynomial of $Y^*Y$ at zero is zero. Since $Y^*Y$ is PSD, and its characteristic polynomial
has three real nonnegative roots, it must have a multiple eigenvalue at zero, and hence its rank cannot be
bigger than 1. Therefore both $Y$ and $Z$ are vectors.

Since the real parts of the products $W\overline{W}$ and $\overline{W}W$ are equal, the real parts of $ZZ^*$ and $Y^*Y$
must also be equal. But $Re\, Y^*Y = Re\, Y^T \overline Y$, and hence by Lemma \ref{SU2SO4} there exist unit norm
quaternions $\zeta_z,\zeta_y$ such that $Z\zeta_z = \zeta_y Y^T = (\zeta_y Y)^T$. Since multiplication of $Z$
from the right and $Y$ from the left by unit norm quaternions does not change the matrices $ZZ^*,Y^*Y$, we
can assume without loss of generality that $Z = Y^T$.

Inserting this in (\ref{commrel62}) and combining with (\ref{commrel6}), we get $[W_{\alpha},W_{\beta}] +
Z_{\alpha}Z_{\beta}^T - Z_{\beta}Z_{\alpha}^T = 0$ for all index pairs $\alpha,\beta = r,i,j,k$. Let us
define a linear map ${\cal W}_Z$ from the linear span ${\cal L}$ of the matrices $W_r,W_i,W_j,W_k$ to
${\mathbb R}^3$ by ${\cal W}_Z(W_{\alpha}) = Z_{\alpha}$, $\alpha = r,i,j,k$. By bilinearity we then have the
relation $[U,V] + {\cal W}_Z(U){\cal W}_Z(V)^T - {\cal W}_Z(V){\cal W}_Z(U)^T = 0$ for all matrices $U,V \in {\cal L}$. Note that ${\cal W}_Z$
has a 1-dimensional kernel. Let $U \in {\cal L}$ be a generator of this kernel. We have ${\cal W}_Z(U) = 0$
and hence $[U,V] = 0$ for all $V \in {\cal L}$. Therefore $U$ commutes with all matrices in the 5-dimensional
subspace of ${\cal S}(3)$ which is generated by ${\cal L}$ and $I_3$. But then $U$ has to be proportional to
$I_3$, which contradicts the non-reducibility of ${\bf B}$. Hence this case cannot occur.

\medskip

Similarly to case 1,\ we can treat the case when $Y^*Y$ has rank at most 1.

\medskip

2. Both $Y^*Y$ and $ZZ^*$ are of rank 2.

Let $v_z,v_y$ be generators of the kernels of $ZZ^*,Y^*Y$,
i.e.\ vectors that satisfy
 $v_z^*ZZ^* = 0$, $Y^*Yv_y = 0$.
Consider the equation $v_z^*(WM - MW)v_y = 0$, where $M$ is a real symmetric matrix. One sees immediately
that $M = I_3$ satisfies this equation. Furthermore, by (\ref{commrel6}), (\ref{commrel62}) and $Re(ZZ^*) =
Re(Y^*Y)$ we have
 \begin{eqnarray*}
 2v_z^*(WW_r - W_rW)v_y &=& 2v_z^*(i[W_i,W_r] + j[W_j,W_r] +
k[W_k,W_r])v_y =
 v_z^*(-ZZ^* + Y^*Y)v_y = 0, \\
 2v_z^*(WW_i - W_iW)v_y &=& 2v_z^*(-[W_i,W_r] + j[W_j,W_i] -
k[W_i,W_k])v_y =
 v_z^*(-ZZ^*i + iY^*Y)v_y = 0, \\
 2v_z^*(WW_j - W_jW)v_y &=& 2v_z^*(-[W_j,W_r] - i[W_j,W_i] +
k[W_k,W_j])v_y =
 v_z^*(-ZZ^*j + jY^*Y)v_y = 0, \\
 2v_z^*(WW_k - W_kW)v_y &=& 2v_z^*(-[W_k,W_r] + i[W_i,W_k] -
j[W_k,W_j])v_y =
 v_z^*(-ZZ^*k + kY^*Y)v_y = 0.
 \end{eqnarray*}
Hence there exists a 5-dimensional linear subspace in ${\cal S}(3)$ that satisfies $v_z^*(WM - MW)v_y = 0$.
This subspace contains the identity matrix and thus intersects the PSD cone $S_+(3)$. But then it contains a
rank 1 PSD matrix $\xi\xi^T$ by Dines theorem \cite{Dines43}. Here $\xi \in {\mathbb R}^3$ is a nonzero
vector.

\medskip

 2.1. $v_z^*\xi = \xi^Tv_y = 0$.

Then there exist nonzero vectors $\zeta_y,\zeta_z$ such that $Z\zeta_z = \xi$ and $Y^*\zeta_y = \xi$. Let
$\alpha = \min(|\zeta_z|^{-2},|\zeta_y|^{-2})$. Then we have $\tilde Z = ZZ^* - \alpha\xi\xi^T = Z(I_2 -
\alpha\zeta_z\zeta_z^*)Z^* \succeq 0$, $\tilde Y = Y^*Y - \alpha\xi\xi^T = Y^*(I_2 - \alpha\zeta_y\zeta_y^*)Y
\succeq 0$. It follows that
\[ {\bf B} - \left( \begin{array}{cc} 0 & 0 \\ 0 &
\alpha\xi\xi^T \end{array} \right) = \left( \begin{array}{cc} I_3 & \overline{W} \\ W & W\overline{W} + \tilde Z
\end{array} \right) \succeq 0, \
\overline{{\bf B}} - \left( \begin{array}{cc} 0 & 0 \\ 0 & \alpha\xi\xi^T \end{array} \right) = \left(
\begin{array}{cc} I & W \\ \overline{W} & \overline{W}W + \tilde Y
\end{array} \right) \succeq 0.
\]
The matrix which is subtracted from ${\bf B}$ and $\overline{{\bf B}}$ is clearly $K_6 \otimes
S_+(3)$-separable and hence represents an element of $\Sigma_{6,3}$. The difference is again a PPT matrix and
hence represents an element in $\Gamma_{6,3}$.

We have represented the element ${\bf B}$ as convex combination of a nonzero extremal element of $\Sigma_{6,3}$ and
another element of $\Gamma_{6,3}$. Since ${\bf B}$ is extremal in $\Gamma_{6,3}$, it must be contained in
$\Sigma_{6,3}$.

\medskip

2.2. $v_z^*\xi \not= 0$ or $\xi^Tv_y \not= 0$.

 We have $v_z^*(W\xi\xi^T - \xi\xi^TW)v_y = 0$.
 Without loss of generality, let $v_z^*\xi \not= 0$. Define $q =
(v_z^*\xi)^{-1}v_z^*W\xi$.
 (In the case $\xi^Tv_y \not= 0$ we define $q =
\xi^TWv_y(\xi^Tv_y)^{-1}$ and proceed
 similarly.)
 Inserting $(v_z^*\xi)q$ for $v_z^*W\xi$, we obtain $\xi^T(W - qI_3)v_y =
0$. Hence we find
 a vector $\zeta_y \in {\mathbb H}^2$ such that $\xi^T(W - qI_3) +
\zeta_y^*Y = 0$, or
 $(\overline{W} - \bar qI_3)\xi + Y^*\zeta_y = 0$. Moreover, $\xi$ commutes with
quaternion scalars and we have
 $v_z^*(W - qI_3)\xi = v_z^*W\xi - v_z^*(v_z^*\xi)^{-1}(v_z^*W\xi)\xi = 0$.
 Hence we find a vector $\zeta_z \in {\mathbb H}^2$ such that $(W -
qI_3)\xi + Z\zeta_z = 0$.
 Let $\alpha = \min((|\xi|^2 + |\zeta_z|^2)^{-1},(|\xi|^2 +
|\zeta_y|^2)^{-1})$.
 Then we have
 \[ {\bf B} - \alpha \left( \begin{array}{c} \xi \\ q\xi
\end{array} \right)
 \left( \begin{array}{c} \xi \\ q\xi \end{array} \right)^* = \left(
\begin{array}{cc} I_3 & 0 \\ W & Z
 \end{array} \right) \left( I_5 - \alpha \left( \begin{array}{c} \xi \\
\zeta_z \end{array} \right) \left( \begin{array}{c} \xi \\ \zeta_z \end{array} \right)^* \right) \left(
\begin{array}{cc} I_3 & 0 \\ W & Z \end{array} \right)^* \succeq 0,
\]
 \[ \overline{{\bf B}} - \alpha \left( \begin{array}{c}
\xi \\ \bar q\xi \end{array} \right)
 \left( \begin{array}{c} \xi \\ \bar q\xi \end{array} \right)^* =
 \left( \begin{array}{cc} I_3 & 0 \\ \overline{W} & Y^* \end{array} \right)
 \left( I_5 - \alpha \left( \begin{array}{c} \xi \\ \zeta_y \end{array}
\right) \left( \begin{array}{c} \xi \\ \zeta_y \end{array} \right)^* \right) \left( \begin{array}{cc} I_3 & 0
\\ \overline{W} & Y^* \end{array} \right)^* \succeq 0.
\]

Thus, as in the previous case, we have represented ${\bf B}$ as
convex combination of an extremal element of $\Sigma_{6,3}$ and a nonzero
element of $\Gamma_{6,3}$. By the extremality of ${\bf B}$ in
$\Gamma_{6,3}$ we obtain ${\bf B} \in \Sigma_{6,3}$.

\medskip

We have proven that any extremal ray of $\Gamma_{6,3}$ is in
$\Sigma_{6,3}$. We obtain the following theorem.

{\theorem \label{63} $\Gamma_{6,3} = \Sigma_{6,3}$. \qed }

Thus a matrix in ${\cal Q}(2) \otimes {\cal S}(3)$ is $Q_+(2) \otimes S_+(3)$-separable if and only if it satisfies the PPT condition.
However, positivity alone is not sufficient for separability, as the following example shows. The matrix
\[ \left( \begin{array}{cccccc} 2 & 0 & 0 & -1 & -i & k \\ 0 & 5 & 0 & -i & 0 & -j \\ 0 & 0 & 2 & k & -j & 1 \\
-1 & i & -k & 4 & 0 & 0 \\ i & 0 & j & 0 & 1 & 0 \\ -k & j & 1 & 0 & 0 & 4 \end{array} \right) \in {\cal Q}(2) \otimes {\cal S}(3)
\]
is PSD, but its transpose and hence its partial transpose is not.

{\corollary \label{PSDPPT} Let $m \geq 2$, $n \geq 3$. Then the cone of PSD matrices in ${\cal Q}(m) \otimes {\cal S}(n)$
is strictly larger than the cone of PPT matrices. \qed }

\section{Counterexample against $\Gamma_{4,4} = \Sigma_{4,4}$}

In this section we present an element of $\Gamma_{4,4}$ that is not in
$\Sigma_{4,4}$.

Recall that $\Gamma_{4,n}$ is the cone of
blockwise symmetric complex hermitian $2n \times 2n$ PSD matrices,
while $\Sigma_{4,n}$ is the cone of $H_+(2) \otimes
S_+(n)$-separable matrices.

Consider the matrix $B = VV^* \in H_+(8)$ with
\[ V = \left( \begin{array}{cc} I_4 & 0 \\ W & Z \end{array} \right),
\qquad W = \left(
\begin{array}{cccc}
-1-i & 2i & i & 0 \\
2i & 1+i & -1+2i & -2i \\
i & -1+2i & 1+i & 2i \\
0 & -2i & 2i & 2-i \end{array} \right) = W^T, \qquad Z = \left(
\begin{array}{c} 1 \\ 3i \\ 0 \\ 0 \end{array} \right).
\]
It is not hard to check that
\[ WW^* + ZZ^* = \left( \begin{array}{cccc}
8 & 2 & 4 & -2 \\ 2 & 24 & 0 & 4 \\ 4 & 0 & 12 & -4 \\ -2 & 4 & -4
& 13
\end{array} \right)
\]
is symmetric and hence $B$ is blockwise symmetric.

Let us show that $B$ is not $H_+(2) \otimes S_+(4)$-separable. It suffices to show that there does not exist
a separable rank 1 matrix $\Xi$ such that $B - \Xi \succeq 0$, i.e.\ there is no nonzero vector $v \in
{\mathbb C}^5$ such that $Vv$ is representable as a single tensor product in ${\mathbb C}^2 \otimes {\mathbb
R}^4$.

Let us assume such a vector $v$ exists. If its first 4 entries are
zero, then $Vv$ is not separable because $Z$ is not a multiple of
a real vector. Hence the first 4 entries of $v$ form a multiple of
a nonzero real vector $v_1$ and the last entry is a complex number
$v_2 = v_2^r + iv_2^i$. Without restriction of generality we can
assume that the first 4 entries of $v$ actually form $v_1$. We
have $Wv_1 + Zv_2 = \lambda v_1$ for some complex number $\lambda
= \lambda^r + i\lambda^i$. Resolving with respect to $v_1,v_2^r,v_2^i$
and separating the real and the imaginary part, we obtain the
homogeneous linear system
\[ \left( \begin{array}{cccc} Re W - \lambda^r I_4 & Re Z & -Im Z \\ Im
W - \lambda^i I_4 & Im Z & Re Z \end{array} \right) \left(
\begin{array}{c} v_1 \\ v_2^r \\ v_2^i \end{array} \right) = 0.
\]
Hence the $8\times 6$ coefficient matrix on the left-hand side
must be rank deficient. The determinant of the matrix formed of
the last six rows is proportional to $\lambda^r(2-\lambda^r)$,
hence $\lambda^r = 0$ or $\lambda^r = 2$. In the first case adding
and in the second case subtracting rows 2 and 3 leads to $v_2^i =
0$. This yields the reduced systems
\[ \left( \begin{array}{ccccc}
-1 & 0 & 0 & 0 & 1 \\
0 & -1 & 1 & 0 & 0 \\
0 & 0 & 0 & 2 & 0 \\
-1-\lambda^i & 2 & 1 & 0 & 0 \\
2 & 1-\lambda^i & 2 & -2 & 3 \\
1 & 2 & 1-\lambda^i & 2 & 0 \\
0 & -2 & 2 & -1-\lambda^i & 0
\end{array} \right) \left( \begin{array}{c} v_1 \\ v_2^r \end{array}
\right) = 0, \]
\[ \left( \begin{array}{ccccc}
-3 & 0 & 0 & 0 & 1 \\
0 & -1 & -1 & 0 & 0 \\
-1-\lambda^i & 2 & 1 & 0 & 0 \\
2 & 1-\lambda^i & 2 & -2 & 3 \\
1 & 2 & 1-\lambda^i & 2 & 0 \\
0 & -2 & 2 & -1-\lambda^i & 0
\end{array} \right) \left( \begin{array}{c} v_1 \\ v_2^r \end{array}
\right) = 0,
\]
which are easily seen to have no nontrivial solution for any $\lambda^i$.

Hence $B$ is not separable and we obtain the following theorem.

{\theorem \label{th44} $\Gamma_{4,4} \not= \Sigma_{4,4}$. \qed }

Theorem \ref{th44} and Lemma \ref{reduceimpl} yield the relation $\Gamma_{4,n} \not= \Sigma_{4,n}$ for any $n \geq 4$ and thus the following result.

{\corollary \label{cor44} The positivity condition in the space ${\cal H}(2) \otimes {\cal S}(n)$ is not sufficient for $H_+(2) \otimes S_+(n)$-separability for $n \geq 4$. \qed }

\smallskip

{\it Remark:} Woronowicz already presented in \cite{Woronowicz} a PPT matrix in $H_+(8)$ which is not $H_+(2)
\otimes H_+(4)$-separable. However, this matrix was not blockwise symmetric and hence not a counterexample
against the equality $\Gamma_{4,4} = \Sigma_{4,4}$.

Theorem \ref{th44} and Lemma \ref{reduceimpl} yield the relation $\Gamma_{6,n} \not= \Sigma_{6,n}$ for any $n \geq 4$ and thus the following result.

{\corollary \label{cor44q} The PPT condition in the space ${\cal Q}(2) \otimes {\cal S}(n)$ is not sufficient for $Q_+(2) \otimes S_+(n)$-separability for $n \geq 4$. \qed }

\section{Positivity implies $S_+(2) \otimes Q_+(n)$-separability}

In this section we investigate the spaces ${\cal S}(2) \otimes {\cal Q}(n)$ for $n \geq 1$.
The dimensionality argument Lemma \ref{dimcount} will be essentially sufficient to prove that the positivity of a matrix is equivalent to
$S_+(2) \otimes Q_+(n)$-separability.

Let $P_n$ be the intersection of the space ${\cal S}(2) \otimes {\cal Q}(n)$ with the cone $Q_+(2n)$ and let $\Sigma_n$ be the $S_+(2) \otimes Q_+(n)$-separable cone.
By Proposition \ref{pospptnec} we have $\Sigma_n \subset P_n$ for any $n \geq 1$ and trivially $P_1 = \Sigma_1$.
Any matrix ${\bf B} \in {\cal S}(2) \otimes {\cal Q}(n)$ can be partitioned as
\begin{equation} \label{quatpart}
{\bf B} = \left( \begin{array}{cc} B_{11} & B_{12} \\ B_{12} & B_{22} \end{array} \right),
\end{equation}
where $B_{11},B_{12},B_{22} \in {\cal Q}(n)$.

{\lemma \label{dimquat} Let ${\bf B}$ generate an extreme ray of $P_n$. Then the rank of ${\bf B}$ is at most $n$. }

\begin{proof}
We apply Lemma \ref{dimcount}. With the notations of this lemma, define $K$ to be the cone $Q_+(2n)$ and $L$ the space ${\cal S}(2) \otimes {\cal Q}(n)$.
Since ${\bf B}$ generates an extreme ray of $P_n$, the dimension of its face in $Q_+(2n)$ can be at most
$(8n^2-2n) - (6n^2-3n) + 1 = 2n^2 + n + 2$. Now suppose that ${\bf B}$ has rank at least $n+1$. Then the dimension of its face in $Q_+(2n)$ equals at least
the dimension of ${\cal Q}(n+1)$, which is $2(n+1)^2 - (n+1) = 2n^2 + 3n + 1 > 2n^2 + n + 2$. This contradiction completes the proof.
\end{proof}

Let us now define a group homomorphism $H_n^{\mathbb H}: GL_n({\mathbb H}) \to \Aut(Q_+(n))$. It shall assign to any regular quaternionic $n \times n$ matrix $S$
the automorphism $A \mapsto SAS^*$. Note that for any $S \in GL_n({\mathbb H})$, the map $\id_{{\cal S}(2)} \otimes H_n^{\mathbb H}(S)$ is both in
$\Aut(P_n)$ and in $\Aut(\Sigma_n)$.

{\lemma Let $A,B \in Q_+(n)$. Then there exists a matrix $S \in GL_n({\mathbb H})$ such that both $SAS^*$ and $SBS^*$ are diagonal. }

\begin{proof}
Assume the conditions of the lemma. Let $l \leq n$ be the rank of $A+B$. Since $A+B \succeq 0$, there exists a regular matrix $S_1$ such that
$S_1(A+B)S_1^*$ is equal to $\diag(I_l,0_{n-l})$.
Now note that $S_1AS_1^* \preceq S_1(A+B)S_1^*$, $S_1BS_1^* \preceq S_1(A+B)S_1^*$ and hence $S_1AS_1^* = \diag(A_1,0_{n-l})$,
$S_1BS_1^* = \diag(B_1,0_{n-l})$ for some matrices $A_1,B_1 \in Q_+(l)$. Let now $U$ be a hyperunitary $l \times l$ matrix such that $U(A_1-B_1)U^*$
is diagonal (see Proposition \ref{spectfac} in the appendix). Then one easily sees that the matrix $S = \diag(U,I_{n-l})S_1$ satisfies the assertion of the lemma.
\end{proof}

{\corollary Let ${\bf B} \in P_n$. Then there exists $S \in GL_n({\mathbb H})$ such that both the upper left and the lower right $n \times n$ block of the matrix
${\bf B}' = (\id_{{\cal S}(2)} \otimes H_n^{\mathbb H}(S))({\bf B})$ are diagonal. \qed }

{\lemma Let ${\bf B}$ generate an extreme ray of $P_n$ and be partitioned as in (\ref{quatpart}). Then there exists $S \in GL_n({\mathbb H})$ such that the map $H_n^{\mathbb H}(S)$ diagonalizes
all three blocks $B_{11},B_{12},B_{22}$. }

\begin{proof}
We prove the lemma by induction. For $n = 1$ the assertion of the lemma holds trivially. Let us now assume that it holds for $n-1$.

Assume the conditions of the lemma. By the preceding corollary we can assume without restriction of generality that $B_{11}$ and $B_{22}$ are already diagonal.

Suppose that the rank of $B_{22}$ is strictly smaller than $n$. Then at least one diagonal element of $B_{22}$, say the last one, is zero. Since ${\bf B} \succeq 0$,
the last row and the last column of ${\bf B}$ are zero. But then also the last row and the last column of $B_{12}$ are zero. Hence all three matrices $B_{11},B_{12},B_{22}$
are block-diagonal, with an upper left block of size $(n-1) \times (n-1)$ and a lower right block of size $1 \times 1$. By assumption of the induction there exists
$S' \in GL_{n-1}({\mathbb H})$ such that the map $H_{n-1}^{\mathbb H}(S)$ diagonalizes all the upper left blocks. It then follows that the matrix
$S = \diag(S',1)$ satisfies the assertion of the lemma.

Let us now assume that $B_{22}$ has full rank. Then there exists $S_1 \in GL_n({\mathbb H})$ such that $S_1B_{22}S_1^* = I_n$. Define 
${\bf B}' = (\id_{{\cal S}(2)} \otimes H_n^{\mathbb H}(S_1))({\bf B})$ and $B_{11}' = S_1B_{11}S_1^*$, $B_{12}' = S_1B_{12}S_1^*$.
Let now $U$ be a hyperunitary $n \times n$ matrix such that $UB_{12}'U^*$ is diagonal.
By Lemma \ref{dimquat} the rank of ${\bf B}'$ equals $n$. This implies $B_{11}' = (B_{12}')^2$. It follows that $UB_{11}'U^* = (UB_{12}'U^*)^2$ is also diagonal.
Thus the matrix $S = US_1$ satisfies the assertion of the lemma.
\end{proof}

{\corollary Let ${\bf B}$ generate an extreme ray of $P_n$. Then ${\bf B} \in \Sigma_n$. \qed }

This yields the following result.

{\theorem \label{th2q} $P_n = \Sigma_n$ for any $n \geq 1$. \qed }

\section{Conclusions}

Let us summarize our results.
As mentioned in the introduction, the cone of $S_+(m) \otimes S_+(n)$-separable matrices in ${\cal S}(m) \otimes {\cal S}(n)$ coincides with the cone of PSD matrices
in this space if and only if $\min(n,m) \leq 2$. This is a consequence of the results in \cite{Terpstra}. 
The cone of $H_+(m) \otimes H_+(n)$-separable matrices in ${\cal H}(mn) = {\cal H}(m) \otimes {\cal H}(n)$ coincides with the cone of PPT matrices
in ${\cal H}(mn)$ if and only if $\min(n,m) = 1$ or $m+n \leq 5$. This is a consequence of the results in \cite{Woronowicz}. We illustrate this in the following tables.

\bigskip

\begin{minipage}[c]{5cm}
    \begin{tabular}{c||c|c|c}
    m\verb1\1n & 2 & 3 & $\geq 4$ \\
    \hline\hline
    2 & PSD & PSD & PSD \\
    \hline
    3 & PSD & N & N \\
    \hline
    $\geq 4$ & PSD & N & N
    \end{tabular}
    
    \medskip
    
    Semidefinite descriptions of the $S_+(m) \otimes S_+(n)$-separable cone
\end{minipage} \hfill
\begin{minipage}[c]{2.5cm}
\end{minipage} \hfill
\begin{minipage}[c]{5cm}
    \begin{tabular}{c||c|c|c}
    m\verb1\1n & 2 & 3 & $\geq 4$ \\
    \hline\hline
    2 & PPT & PPT & N \\
    \hline
    3 & PPT & N & N \\
    \hline
    $\geq 4$ & N & N & N
    \end{tabular}
    
    \medskip
    
    Semidefinite descriptions of the $H_+(m) \otimes H_+(n)$-separable cone
\end{minipage}

\bigskip

Here "PSD" indicates that the cone of separable elements in the space corresponding to the pair $(m,n)$ equals the cone of positive semidefinite matrices,
and "PPT" indicates that the cone of separable elements equals the cone of PPT matrices. An "N" indicates that the separable cone is described by neither
the PSD cone nor the PPT cone.

\medskip

Theorems \ref{thH2S}, \ref{C43}, Corollary \ref{cor44} and Lemma \ref{33conseq} lead to the following result.

{\theorem The cone of $H_+(m) \otimes S_+(n)$-separable matrices in ${\cal H}(m) \otimes {\cal S}(n)$ coincides with the cone of PSD matrices
in this space if and only if $n \leq 2$ or $m = 1$ or $m+n \leq 5$. \qed }

\medskip

Theorem \ref{th2q}, Lemma \ref{33conseq}, and Corollary \ref{PSDPPT} lead to the following result.

{\theorem The cone of $Q_+(m) \otimes S_+(n)$-separable matrices in ${\cal Q}(m) \otimes {\cal S}(n)$ coincides with the cone of PSD matrices
in this space if and only if $n \leq 2$ or $m = 1$. \qed }

\medskip

Theorem \ref{63}, Corollary \ref{cor44q} and Lemma \ref{33conseq} lead to the following result.

{\theorem The cone of $Q_+(m) \otimes S_+(n)$-separable matrices in ${\cal Q}(m) \otimes {\cal S}(n)$ coincides with the cone of PPT matrices
in this space if and only if $m \leq 2$ and $m+n \leq 5$. \qed }

We can summarize these results in the following tables.

\bigskip

\begin{minipage}[c]{5cm}
    \begin{tabular}{c||c|c|c}
    m\verb1\1n & 2 & 3 & $\geq 4$ \\
    \hline\hline
    2 & PSD & PSD & N \\
    \hline
    3 & PSD & N & N \\
    \hline
    $\geq 4$ & PSD & N & N
    \end{tabular}
    
    \medskip
    
    Semidefinite descriptions of the $H_+(m) \otimes S_+(n)$-separable cone
\end{minipage} \hfill
\begin{minipage}[c]{2.5cm}
\end{minipage} \hfill
\begin{minipage}[c]{5cm}
    \begin{tabular}{c||c|c|c}
    m\verb1\1n & 2 & 3 & $\geq 4$ \\
    \hline\hline
    2 & PSD & PPT & N \\
    \hline
    3 & PSD & N & N \\
    \hline
    $\geq 4$ & PSD & N & N
    \end{tabular}
    
    \medskip
    
    Semidefinite descriptions of the $Q_+(m) \otimes S_+(n)$-separable cone
\end{minipage}

\bigskip

In all four tables, for $\min(n,m) = 1$ the cone of separable matrices is trivially equal to the cone of PSD matrices.

The preceding two theorems yield also the following nontrivial result.

{\corollary A matrix in ${\cal Q}(2) \otimes {\cal S}(2)$ is PSD if and only if its transpose is PSD. \qed }

Apart from the trivial cases ${\cal Q}(2) = {\cal Q}(2) \otimes {\cal S}(1)$ and ${\cal S}(n) = {\cal Q}(1) \otimes {\cal S}(n)$ the
space ${\cal Q}(2) \otimes {\cal S}(2)$ is thus the only tensor product space ${\cal Q}(m) \otimes {\cal S}(n)$, $m,n \in {\mathbb N}_+$,
where the positivity property is invariant with respect to transposition (cf.\ Corollaries \ref{posQ2} and \ref{posQn} in the appendix).

\bibliography{convexity,matrices,quaternions,quantinf,geometry}

\begin{thebibliography}{10}

\bibitem{KoecherBraun}
Hel Braun and Max Koecher.
\newblock {\em Jordan-Algebren}, volume 128 of {\em A Series of Comprehensive
  Studies in Mathematics}.
\newblock Springer, Berlin, New York, 1966.

\bibitem{Brenner51}
J.L. Brenner.
\newblock Matrices of quaternions.
\newblock {\em Pacific J. Math.}, 1:329--335, 1951.

\bibitem{Choi74}
M.-D. Choi.
\newblock A {S}chwarz inequality for positive linear maps on {C}-algebras.
\newblock {\em Illinois J. Math.}, 18(4):565--574, 1974.

\bibitem{Dines43}
Lloyd~L. Dines.
\newblock On linear combinations of quadratic forms.
\newblock {\em Bull. Amer. Math. Soc.}, 49:388--393, 1943.

\bibitem{Gurvits0302102}
Leonid Gurvits and Howard Barnum.
\newblock Separable balls around the maximally mixed multipartite quantum
  states.
\newblock {\em Phys. Rev.}, 68:042312, 2003.

\bibitem{Horodeckis96}
Ryszard Horodecki, Pawel Horodecki, and Michal Horodecki.
\newblock Separability of mixed states: necessary and sufficient conditions.
\newblock {\em Physics Letters}, 223(1):1--8, 1996.

\bibitem{Lee49}
H.C. Lee.
\newblock Eigenvalues and canonical forms of matrices with quaternion
  coefficients.
\newblock {\em Proceedings of the Royal Irish Academy, Sect. A}, 52:253--260,
  1949.

\bibitem{Peres96}
A.~Peres.
\newblock Separability criterion for density matrices.
\newblock {\em Phys. Rev. Lett.}, 77:1413--1415, 1996.

\bibitem{Stormer63}
Erling St{\o}rmer.
\newblock Positive linear maps of operator algebras.
\newblock {\em Acta Mathematica}, 110:233--278, 1963.

\bibitem{Terpstra}
F.J. Terpstra.
\newblock {Die Darstellung biquadratischer Formen als Summen von Quadraten mit
  Anwendung auf die Variationsrechnung}.
\newblock {\em Mathematische Annalen}, 116:166--180, 1939.

\bibitem{Vinberg63}
E.~B. Vinberg.
\newblock The theory of convex homogeneous cones.
\newblock {\em Transactions of Moscow Mathematical Society}, 12:340--403, 1963.

\bibitem{Warner}
Frank Warner.
\newblock {\em Foundations of Differentiable Manifolds and {L}ie Groups}.
\newblock Springer, New York, Berlin, Heidelberg, Tokyo, 1983.

\bibitem{Wiegmann55}
N.A. Wiegmann.
\newblock Some theorems on matrices with real quaternion elements.
\newblock {\em Can. J. Math.}, 7:191--201, 1955.

\bibitem{Wolf36}
L.A. Wolf.
\newblock Similarity of matrices in which the elements are real quaternions.
\newblock {\em Bulletin of the American Mathematical Society}, 42:737--743,
  1936.

\bibitem{Woronowicz}
S.L. Woronowicz.
\newblock Positive maps of low dimensional matrix algebras.
\newblock {\em Reports on Mathematical Physics}, 10:165--183, 1976.

\bibitem{Yakubovitch70}
V.A. Yakubovitch.
\newblock Factorization of symmetric matrix polynomials.
\newblock {\em Doklady Akademii Nauk SSSR}, 194(3):1261--1264, 1970.

\bibitem{Zhang}
Fuzhen Zhang.
\newblock Quaternions and matrices of quaternions.
\newblock {\em Linear algebra and its applications}, 251:21--57, 1997.

\end{thebibliography}
\bibliographystyle{plain}

\appendix

\section{Quaternions}

In this section we provide some basic facts about quaternions and
quaternionic matrices, which will be used in the paper. 
Most concepts from the theory of real and
complex hermitian matrices carry over to the quaternionic case. For a detailed
treatment of the properties
of quaternionic matrices see \cite{Brenner51},\cite{Zhang}.

\medskip

We denote the algebra of quaternions by ${\mathbb H}$, with the usual
product rules
\begin{equation} \label{quatgen}
i^2 = j^2 = k^2 = ijk = -1
\end{equation}
for the generators. This algebra is not
commutative. However, any non-zero quaternion $q = q_r + iq_i + jq_j +
kq_k$ has a two-sided inverse
$q^{-1} = \bar q/|q|^2$, where $\bar q = q_r - iq_i - jq_j - kq_k$ is
the quaternionic conjugate of $q$
and $|q| = \sqrt{q_r^2+q_i^2+q_j^2+q_k^2}$ is the norm of $q$. The real
numbers $q_r,q_i,q_j,q_k$ are the
components of $q$.
The norm has the properties $|q| = |\bar q|$ and
$|pq| = |p||q|$. The conjugate of a product is given by $\overline{pq} =
\bar q \bar p$. Denote by $Re\,q$ the real
part $(q+\bar q)/2$ of $q$.

Let ${\cal Q}_{ij}: {\mathbb H} \to {\mathbb H}$ be the map defined by
${\cal Q}_{ij}: q = q_r + iq_i + jq_j + kq_k \mapsto q_r - iq_i - jq_j +
kq_k$. It is not
hard to see that ${\cal Q}_{ij}$ is an automorphism of the quaternion
algebra, since it leaves relations
(\ref{quatgen}) invariant. Likewise
we can define the automorphisms ${\cal Q}_{ik}$ and ${\cal Q}_{jk}$,
which switch the
signs of the corresponding components. There exist also other
automorphisms, e.g.\
the cyclic permutation of the imaginary units $i,j,k$.

\medskip

Let us now consider vectors and matrices of quaternions.

Denote the space of column vectors with $n$ quaternionic entries
by ${\mathbb H}^n$ and the space of quaternionic matrices of size $m \times n$ by
${\mathbb H}^{m \times n}$.
The norm $||v||$ of a quaternionic vector $v$ can be defined as $\sqrt{v^*v}$.
Here $v^*$ is the quaternionic conjugate transpose of $v$.
This norm equals the Euclidean norm in the associated real vector space
${\mathbb R}^{4n}$ of the components
of the elements of $v$.

The conjugate $\overline A$ of a quaternionic matrix $A$ is defined entrywise.
Let $A^*$ denote the conjugate transpose of $A$.
Then for matrix products we obviously have $(AB)^* = B^*A^*$.

The quaternionic matrix algebras have a representation as real matrix
algebras.
Namely, if we write quaternionic matrices $A$ as $A = A_r + iA_i + jA_j
+ kA_k$,
where $A_r,A_i,A_j,A_k$ are real matrices, then the isomorphism is given by
\begin{equation} \label{quatiso}
{\cal I}_q:\quad A \mapsto \left( \begin{array}{cccc} A_r & -A_i & -A_j
& -A_k
\\ A_i & A_r & -A_k & A_j \\ A_j & A_k & A_r & -A_i \\ A_k & -A_j
& A_i & A_r \end{array} \right).
\end{equation}
The conjugate transpose of a quaternionic matrix corresponds
to the transpose of its real counterpart, ${\cal I}_q(A^*) = ({\cal
I}_q(A))^T$.

The trace of matrix products is not invariant under cyclic permutations,
even for scalars.
Nevertheless, we have the following weaker assertion.

{\prop For quaternionic matrices $A,B$ of
appropriate sizes the relation $Re(tr(AB)) = Re(tr(BA))$ holds. }

\begin{proof}
For a square matrix $A$ we have $tr {\cal I}_q(A) = 4 Re(tr A)$. Hence
$Re(tr (AB)) = \frac{1}{4}tr({\cal I}_q(A) {\cal I}_q(B)) =
\frac{1}{4}tr({\cal I}_q(B) {\cal I}_q(A)) =
Re(tr (BA))$.
\end{proof}

Let $A$ be an $n \times n$ square quaternionic matrix.
If there does not exist a nonzero quaternionic vector $v$ such
that $Av = 0$, then there exists an inverse quaternionic matrix
$A^{-1}$ with $AA^{-1} = A^{-1}A = I_n$ \cite{Brenner51}. We will call
such matrices {\sl regular}.
A regular matrix $A$ can be viewed as an ${\mathbb R}$-linear
automorphism of ${\mathbb H}^n$, which maps
a vector $v$ to the vector $Av$.
If there exists $v \not= 0$ such that $Av = 0$ then we call $A$ {\sl
singular}.
Obviously ${\cal I}_q(A^{-1}) = ({\cal I}_q(A))^{-1}$, and the inverse
$A^{-1}$ of a regular
matrix is unique. The inverse of a product is given by $(AB)^{-1} =
B^{-1}A^{-1}$, the
inverse of the conjugate transpose by $(A^*)^{-1} = (A^{-1})^*$, which
is checked easily by applying
the isomorphism ${\cal I}_q$.

If for a square quaternionic matrix $A$ and a quaternionic vector $v$ there
exists a quaternion $\lambda$ such that $Av = v\lambda$, then we call $v$
{\sl right eigenvector} of $A$ and $\lambda$ the corresponding {\sl
right eigenvalue}.
An $n\times n$ quaternionic matrix has $n$ right eigenvalues, which are
determined up to
similarity transformations (i.e.\ transformations
$\lambda \mapsto q\lambda q^{-1}$ for $q \in {\mathbb H}$, $q \not= 0$).
The right eigenvalues are preserved under
similarity transformations of the matrix (i.e.\ $A \mapsto SAS^{-1}$ for
regular $S$).
If $A$ is triangular, then the diagonal elements are representatives of
the (left and right) eigenvalues \cite{Lee49}.

If for a square matrix $U$ we have $UU^* = I$, then we call $U$ {\sl
hyperunitary}.
A matrix $U$ is hyperunitary if and only if the matrix ${\cal I}_q(U)$ is
orthogonal, hence the hyperunitary
matrices of size $n \times n$ form a compact group, the {\sl compact symplectic group} $Sp(n)$. 
Hyperunitary matrices can be viewed as
norm-preserving ${\mathbb R}$-linear automorphisms of ${\mathbb H}^n$
(though not every such
automorphism can be represented by a hyperunitary matrix).

Any quaternionic matrix $A$ can be decomposed as $A = UDV$, where $U,V$
are hyperunitary
matrices and $D$ is a diagonal matrix, which has the same size as $A$,
with real nonnegative
entries. The diagonal elements are called {\sl singular values} of $A$
\cite{Lee49},
\cite{Wiegmann55}. For a quaternionic matrix $A$, we call the number of its
positive singular values the {\sl rank} of $A$. Obviously an $n \times
n$ matrix is regular
if and only if it has rank $n$. The rank is invariant under
multiplications by regular
matrices from the left and from the right \cite{Zhang}. It can be shown
that the rank of a quaternionic
matrix is equal to the maximal number of right linearly (over ${\mathbb
H}$) independent columns
and left linearly independent rows of the matrix \cite{Wolf36}.

We call a square quaternionic matrix $A$ {\sl hermitian} if $A = A^*$.
The hermitian matrices of size $n \times n$ form a $(2n^2 -
n)$-dimensional vector space ${\cal Q}(n)$ over the
reals. The real part $A_r$ of a hermitian quaternionic matrix is symmetric,
while the three imaginary parts $A_i,A_j,A_k$ are skew-symmetric.
There is a scalar product on ${\cal Q}(n)$ given by
\[ \langle A,B \rangle = tr(A_rB_r - A_iB_i - A_jB_j - A_kB_k) =
Re(tr (AB)) = Re(tr (BA)).
\]

{\prop \label{spectfac}
Let $A \in {\cal Q}(n)$. Then there exists $U \in Sp(n)$ such that
$A = UDU^*$, where $D$ is a diagonal matrix with real entries. The columns of $U$ are
right eigenvectors of $A$, while the
diagonal elements of $D$ are the corresponding right eigenvalues. } 

This follows from a well-known
result on factorization of quaternionic matrices \cite{Brenner51},
namely that for any square
quaternionic matrix $A$ there exists a hyperunitary matrix $U$ such that $A =
UTU^*$, where
$T$ is upper triangular. Namely, if $A$ is hermitian, then $T$ must also
be hermitian and
hence diagonal and real. 

The right eigenvalues of a hermitian matrix are real and uniquely determined.
The singular values of a hermitian matrix
are given by the absolute values of its eigenvalues, hence its rank equals
the number of non-zero eigenvalues.
Since the eigenvalues $\lambda_l$ of $A$ are well-defined, its
characteristic
polynomial $p(\lambda) = \prod_{l=1}^n (\lambda - \lambda_l)$ is
well-defined and has
real coefficients. 

{\prop \label{charpolquat} Let $A \in {\cal Q}(n)$. Then the characteristic polynomial of
${\cal I}_q(A)$ is the
fourth power of the characteristic polynomial of $A$. }

\begin{proof}
Consider the matrix ${\cal I}_q(A)$ as an element of
${\cal Q}(4n)$.
Its eigenvalues are real and do not change under conjugation with the
hyperunitary matrix
\[ \frac{1}{2} \left( \begin{array}{cccc} 1 & i & j & k \\ -1 & i & -j &
k \\
-1 & -i & j & k \\ 1 & -i & -j & k \end{array} \right) \otimes I_n.
\]
This conjugation leads to a block-diagonal matrix with blocks
$A,A_r+iA_i-jA_j-kA_k,A_r-iA_i+jA_j-kA_k,A_r-iA_i-jA_j+kA_k$. Hence the
characteristic
polynomial of ${\cal I}_q(A)$ equals the product of the characteristic
polynomials of these
4 blocks. But the last three blocks are the images of $A$ under the
automorphisms
${\cal Q}_{ij},{\cal Q}_{jk},{\cal Q}_{ik}$. Since the real line is
invariant under these
automorphisms, the eigenvalues of these blocks coincide with the
eigenvalues of $A$.
Therefore their characteristic polynomials also coincide with that of
$A$. This completes
the proof.
\end{proof}

It follows that the characteristic polynomial of a hermitian matrix $A$
can be computed by
taking the unique fourth root of the characteristic polynomial of ${\cal
I}_q(A)$ whose leading coefficient equals 1.

A matrix $A \in {\cal Q}(n)$ is called {\sl positive semidefinite} (PSD)
if for any vector $v \in {\mathbb H}^n$
we have $v^* A v \geq 0$ and {\sl positive definite} (PD) if for any non-zero
vector $v \in {\mathbb H}^n$ we have $v^* A v > 0$. Note that $(v^* A v)^* = v^* A^* v$,
hence $v^* A v \in {\mathbb R}$. Clearly $A$ is PSD if and only if all of its
eigenvalues are nonnegative and
PD if and only if all of its eigenvalues are positive.
Let us denote the cone of PSD
matrices in ${\cal Q}(n)$ by $Q_+(n)$.

For any $A \in Q_+(n)$ of rank $l$ there exists
$V \in {\mathbb H}^{n \times l}$ of rank $l$ such that $A = VV^*$.
Namely, let $A = U\left( \begin{array}{cc} D & 0 \\ 0 & 0 \end{array}
\right)U^*$, where $U$
is hyperunitary and $D$ is a real diagonal $l \times l$ matrix with positive
diagonal elements.
Then we can choose $V = U\left( \begin{array}{c} D^{1/2} \\ 0
\end{array} \right)$.
We have the freedom of multiplying the factor $V$ by a hyperunitary matrix
from the right.
On the other hand, if $A = VV^*$, then $A \succeq 0$ and the rank of $A$
equals the rank of $V$. 
Let us define $A^{1/2} = UD^{1/2}U^*$, then $A^{1/2} \in Q_+(n)$ and
$A^{1/2}A^{1/2} = A$. If $A$ is PD, then
$A^{1/2}$ is also PD and its inverse is given by $A^{-1/2} = UD^{-1/2}U^*$.

If for some vector $v \in {\mathbb H}^n$ and a matrix $A \in Q_+(n)$ we have $v^* A
v = 0$, then
$Av = 0$. This is because $v^* A v = (A^{1/2}v)^*(A^{1/2}v) =
|A^{1/2}v|^2 = 0$ yields
$A^{1/2}v = 0$ and therefore $A^{1/2}A^{1/2}v = 0$.

Now let $n = n_1+n_2$. If $A = \left( \begin{array}{cc} A_{11} & A_{12} \\ A_{21} & A_{22}
\end{array} \right) \in Q_+(n)$ is partitioned in four blocks such that
$A_{11} \in {\cal Q}(n_1)$ and $A_{22} \in {\cal Q}(n_2)$
are square matrices, and $v \in {\cal H}^{n_1}$ is such that $A_{11}v
= 0$,
then we have also $A_{21}v = 0$. This can be seen as follows.
If we define $\tilde v = \left( \begin{array}{c} v \\ 0_{n_2 \times 1} \end{array}
\right)$, then
$\tilde v^* A \tilde v = v^* A_{11} v = 0$. Hence $A \tilde v = \left(
\begin{array}{c}
A_{11}v \\ A_{21}v \end{array} \right) = 0$.

Let $\sigma \subset \{1,\dots,n\}$ be a subset of indices, let $A =
VV^*$ be a PSD $n\times n$
matrix. Then the submatrix $\tilde A$ of $A$
that consists of those elements whose row and column numbers are in
$\sigma$ is
also PSD, because it can be represented as product $\tilde V\tilde V^*$,
where $\tilde V$
consists of those elements of $V$ whose row numbers are in $\sigma$.
Moreover, if
$A \succ 0$, then $\tilde A \succ 0$. On the other hand, a
block-diagonal hermitian
matrix is PSD if and only if each block is PSD, and it is PD if and only
if each block is PD.

It is well-known that the cone $Q_+(n)$ is convex and self-dual
\cite{Vinberg63}.
Any automorphism of the quaternion algebra induces an automorphism of
the cone $Q_+(n)$,
because it leaves the real line invariant.

The cone $H_+(n)$ of complex hermitian PSD matrices is an intersection of $Q_+(n)$ with a linear subspace.
Namely, if $A$ is hermitian with zero $j$- and $k$-components, then $A
\in Q_+(n)$
if and only if ${\cal I}_q(A) \in S_+(4n)$ and hence if and only if $A
\in H_+(n)$.

Any PSD quaternionic matrix can be written as sum of
matrices of the form $vv^*$, where $v$ are vectors. These vectors can be
obtained
e.g.\ as columns of the factor $V$ of $A = VV^*$. Thus $Q_+(n)$
is the convex conic hull of the set of rank 1 matrices $\{vv^*\,|\,v \in
{\mathbb H}^n\}$, and any such rank 1 matrix generates an extreme ray of $Q_+(n)$.
The matrix $vv^*$, where $v \in {\mathbb H}^n$ is a nonzero vector, has only one positive
eigenvalue, namely $|v|^2 = v^*v$. Hence $I - vv^* \succeq 0$ for all $v$ in the unit ball.

{\lemma \label{rankiso} Let $S \in {\mathbb H}^{n \times l}$ be of rank
$l$, where $l \leq n$, and let $A \in {\mathbb H}^{l\times l}$. Then
$SAS^* = 0$ is equivalent to $A = 0$. If $A \in {\cal Q}(l)$, then
$SAS^* \succeq 0$ is equivalent to $A \succeq 0$. }

\begin{proof} Clearly $A = 0$ implies $SAS^* = 0$ and $A \succeq 0$
implies $SAS^* \succeq 0$.
Let $S = UDV$ be the singular value decomposition of $S$. Then
$S^*S = V^*D^TDV$. But $D^TD$ is PD, hence $S^*S$ is also PD and invertible.
Then we have $A = (S^*S)^{-1}S^*(SAS^*)S(S^*S)^{-1}$. Hence the converse
implications hold too. \end{proof}

{\lemma Let $S \in {\mathbb H}^{n \times l}$ be of rank $r$. Then there exists
$T \in {\mathbb H}^{(n-r)\times n}$ of rank $n-r$ such that $TS = 0$. Moreover, if for some 
$S' \in {\mathbb H}^{n \times l'}$ we have
$TS' = 0$, then there exists $W \in {\mathbb H}^{l \times l'}$ such that $S' =
SW$. }

\begin{proof}
Let $S = U
\left( \begin{array}{cc} D & 0_{r\times (l-r)} \\ 0_{(n-r)\times r} &
0_{(n-r)\times(l-r)}
\end{array} \right)V$ be the singular
value decomposition of $S$, where $D$ is a positive definite diagonal $r
\times r$ matrix.
Define $T = (0_{(n-r)\times r}\ I_{n-r})U^*$. Then $T$ is of rank $n-r$
and $TS = 0$.
Suppose $S'$ satisfies the assumptions of the lemma. Then the last $n-r$
rows of $U^*S'$ are
zero. Denote the matrix given by the first $r$ rows by $P$. Define
\[ W = V^* \left( \begin{array}{c} D^{-1}P \\ 0_{(l-r)\times l'}
\end{array} \right).
\]
Then we have
\[ SW = U\left( \begin{array}{cc} D & 0_{r\times (l-r)} \\
0_{(n-r)\times r} &
0_{(n-r)\times(l-r)} \end{array} \right)
\left( \begin{array}{c} D^{-1}P \\ 0_{(l-r)\times l'} \end{array} \right) =
U\left( \begin{array}{c} P \\ 0_{(n-r)\times l'} \end{array}
\right) = U(U^*S') = S'.
\quad \qedhere
\]
\end{proof}

Let us investigate the facial structure of $Q_+(n)$.

{\prop \label{quatface} Let $S \in {\mathbb H}^{n \times l}$ be of rank $l$, $l \leq n$.
Then the set $\{SBS^*\,|\,B \succeq 0\}$ is a face of $Q_+(n)$. This
face is isomorphic to $Q_+(l)$. Any face of
$Q_+(n)$ can be expressed in such a way. }

\begin{proof}
Let $A \in Q_+(n)$ be arbitrary. Denote the rank of
$A$ by $l$. Then there exists $S \in {\mathbb H}^{n \times l}$ of rank
$l$ such that $A = SS^*$. By Lemma \ref{rankiso} the map $B \mapsto SBS^*$
is an isomorphism between the space ${\cal Q}(l)$ and the linear
subspace $L = \{SBS^*\,|\,B \in {\cal Q}(l)\} \subset {\cal Q}(n)$, and maps
$Q_+(l)$ to the intersection $L_+ = L \cap Q_+(n)$. Moreover,
$A$ is in the interior of $L_+$.
Let us show that $L_+$ is a face of $Q_+(n)$. 

Let $M_1,M_2 \in Q_+(n)$ such that $\frac{M_1+M_2}{2} = SBS^*$ for some $B \in Q_+(l)$.
We have to show that $M_1,M_2 \in L_+$. Define $C = M_2-M_1 \in {\cal Q}(n)$,
then $SBS^* + \alpha C \in Q_+(n)$ for all $\alpha \in [-1/2,+1/2]$.
By the previous lemma there exists $T \in {\mathbb H}^{(n-l)\times n}$ of
rank $n-l$ such that $TS = 0$. We have $T(SBS^* + \alpha C)T^* =
\alpha TCT^* \succeq 0$ for all $\alpha \in [-1/2,+1/2]$. Hence $TCT^* = 0$
and $T(SBS^* + \alpha C)T^* = 0$ for all $\alpha \in {\mathbb R}$. 
But $SBS^* + \alpha C \succeq 0$ for $\alpha \in [-1/2,+1/2]$, 
therefore we also have $(SBS^* + \alpha C)T^* = \alpha CT^* =
0$ and $TC = 0$. By the previous lemma there exists a matrix $W$
such that $C = SW$. Since $C$ is hermitian, we also have $C =
W^*S^*$. It follows that $W = (S^*S)^{-1}S^*SW =
(S^*S)^{-1}S^*W^*S^*$ and therefore $C = S(S^*S)^{-1}S^*W^*S^* =
SWS(S^*S)^{-1}S^* = S\frac{(S^*S)^{-1}S^*W^* + WS(S^*S)^{-1}}{2}S^*$. 
We have constructed a hermitian matrix $B'$
such that $C = SB'S^* \in L$. Therefore $M_1,M_2 \in L$ and hence $M_1,M_2 \in L_+$.

Let now $S \in {\mathbb H}^{n \times l}$ be given and of rank $l$. Then we define $A = SS^*$ 
and proceed as above.
\end{proof}

Similar results on the facial structure of $S_+(n)$ and $H_+(n)$ can be obtained by the same line of reasoning.

{\lemma \label{quat_schur_compl} Let $n = n_1+n_2$. Let
$A = \left( \begin{array}{cc} A_{11} & A_{12} \\ A_{21} & A_{22}
\end{array} \right) \in {\cal Q}(n)$
be partitioned in four blocks such that $A_{11} \in Q_+(n_1)$ is PD.
Then $A \succeq 0$ if and only if $A_{22} - A_{21}A_{11}^{-1}A_{12} \succeq 0$. 
Moreover, the rank of $A$ is the sum of $\rk(A_{22} - A_{21}A_{11}^{-1}A_{12})$
and $\rk A_{11}$. }

\begin{proof} If $A_{22} - A_{21}A_{11}^{-1}A_{12}$ is PSD of rank $l$, then there
exists $Z \in {\mathbb H}^{n_2 \times l}$ of rank $l$ such that
\[ A = \left( \begin{array}{cc} A_{11}^{1/2} & 0_{n_1 \times l} \\ A_{21}A_{11}^{-1/2}
& Z \end{array}
\right) \left( \begin{array}{cc} A_{11}^{1/2} & 0_{n_1 \times l} \\ A_{21}A_{11}^{-1/2}
& Z \end{array}
\right)^*.
\]
Hence $A$ is PSD.

Let $A \succeq 0$ and let $x_2 \in {\mathbb H}^{n_2}$ be arbitrary. Define $x_1 =
-A_{11}^{-1}A_{12}x_2 \in {\mathbb H}^{n_1}$, then we have
\begin{eqnarray*}
(x_1^*\ x_2^*) A \left( \begin{array}{c} x_1 \\ x_2 \end{array} \right) &=&
x_2^*A_{21}A_{11}^{-1}A_{12}x_2 - x_2^*A_{21}A_{11}^{-1}A_{12}x_2 -
x_2^*A_{21}A_{11}^{-1}A_{12}x_2 + x_2^*A_{22}x_2 \\ &=&
x_2^*(A_{22} - A_{21}A_{11}^{-1}A_{12})x_2 \geq 0,
\end{eqnarray*}
and $A_{22} - A_{21}A_{11}^{-1}A_{12} \succeq 0$.

Now let $A_{11}^{1/2} = UDU^*$ and $Z = U'D'V'$ be the singular value
decompositions of $A_{11}^{1/2}$ and $Z$. Define
\[ P = \left( \begin{array}{cc} U^* & 0 \\ 0 & {U'}^* \end{array} \right)
\left( \begin{array}{cc} A_{11}^{1/2} & 0_{n_1 \times l} \\ A_{21}A_{11}^{-1/2} & Z
\end{array} \right)
\left( \begin{array}{cc} U & 0 \\ 0 & {V'}^* \end{array} \right) =
\left( \begin{array}{cc} D & 0_{n_1 \times l} \\ {U'}^*A_{21}A_{11}^{-1/2}U & D' \end{array} \right) \in {\mathbb H}^{n \times (n_1+l)}
\]
and partition it into matrices $P_1 \in {\mathbb H}^{(n_1+l) \times (n_1+l)}$, $P_2 \in {\mathbb H}^{(n_2-l) \times (n_1+l)}$.
Note that $P$ and $P_1$ are lower triangular matrices, with positive elements on the diagonal. Hence $P_1^*$
is hence invertible.
But then $P_1^*P_1$ and $P^*P = P_1^*P_1 + P_2^*P_2$ are PD. 
Hence the rank of $P$
equals the number $n_1+l$ of its columns, which in turn is equal to
$\rk A_{11} + \rk(A_{22} - A_{21}A_{11}^{-1}A_{12})$. This completes
the proof.
\end{proof}

{\corollary \label{posQ2} A matrix $A = (a_{\alpha\beta}) \in {\cal Q}(2)$ is PSD if
and only if $a_{11} \geq 0, a_{22} \geq 0, a_{11}a_{22} \geq |a_{12}|^2$.
In particular, the cone $Q_+(2)$ is invariant with respect to transposition. \qed }

However, in general $A \succeq 0$ does not yield $\overline A \succeq 0$.
For example, we have
\[
\left( \begin{array}{c} 1 \\ i \\ j \end{array} \right)
\left( \begin{array}{c} 1 \\ i \\ j \end{array} \right)^* = \left(
\begin{array}{ccc} 1 & -i & -j \\ i & 1 & -k \\ j & k & 1
\end{array} \right) \succeq 0,\quad
\left( \begin{array}{c} 1 \\ i \\ j \end{array} \right)^* \left(
\begin{array}{ccc} 1 & i & j \\ -i & 1 & k \\ -j & -k & 1
\end{array} \right) \left( \begin{array}{c} 1 \\ i \\ j \end{array} \right)
= -3 < 0.
\]

{\corollary \label{posQn} The cone $Q_+(n)$ is invariant with respect to transposition
if and only if $n \leq 2$. \qed }

\end{document}